\documentstyle[amssymb,amsfonts,psfig]{amsart}

\newcommand{\rr}{\mbox{{\bf R}}}  % disable Bbb
\newcommand{\norm}[1]{ \left|  #1 \right| }

\newcommand{\Norm}[1]{ \left\| #1 \right\| }

\newcommand{\h}{{h}}

\newcommand{\gamab}{\Gamma_{\alpha \beta}}
\newcommand{\sol}{{S}}

\newcommand{\g}{g}
\newcommand{\m}{{\cal{M}}}

\newcommand{\ddtil}{{\widetilde{\cal D}}}
\newcommand{\xtil}{{\tilde{X}}}
\newcommand{\dd}{{{\cal D}}}

\def\be#1{\begin{equation} \label{#1}}
\def\bs{\begin{split}}
\def\es{\end{split}}
\def\ba{\begin{align}}
\def\bas{\begin{align*}}
\def\ea{\end{align}}
\def\eas{\end{align*}}
\def\R{{\hbox{\bf R}}}

\def\dist{{\hbox{dist}}}

\def\R{{\hbox{\bf R}}}

\def\stil{{\tilde{s}}}

\def\eps{\varepsilon}

\newenvironment{proof}{\noindent {\bf Proof} }{\endprf\par}
\def \endprf{\hfill  {\vrule height6pt width6pt depth0pt}\medskip}
\def\emph#1{{\it #1}}
\def\textbf#1{{\bf #1}}

\theoremstyle{plain}
  \newtheorem{theorem}[subsection]{Theorem}
  
  \newtheorem{proposition}[subsection]{Proposition}
  \newtheorem{lemma}[subsection]{Lemma}
  \newtheorem{corollary}[subsection]{Corollary}

\theoremstyle{remark}

\theoremstyle{definition}
  \newtheorem{definition}[subsection]{Definition}

\include{psfig}

\begin{document}

\title[Wave maps on $\R^{1+1}$]{Local and global well-posedness 
of wave maps on $\R^{1+1}$ for rough data}
\author{Markus Keel, Terence Tao}
\address{Department of Mathematics, UCLA, Los Angeles CA 90095-1555}
\email{keel@@math.ucla.edu, tao@@math.ucla.edu}
\subjclass{35J10,42B25}

\vspace{-0.3in}
\begin{abstract}
We consider wave maps between 
Minkowski space $\rr^{1+1}$ and an analytic manifold. Results  include
global existence for large data in  Sobolev spaces $H^s$ for $s > 3/4$, and
in the scale-invariant norm $L^{1,1}$. We prove local well-posedness in
$H^s$ for $s > 3/4$, and a negative well-posedness result for wave maps
on $\rr^{n+1}$ with data in $H^{n/2}(\rr^{n})$, $n \geq 1$. 
Also included are positive and negative results for scattering. 
 
\end{abstract}

\vspace{-1.0in}

\maketitle

\tableofcontents
\section{Introduction}

Write $(\rr^{n+1}, \g)$ for $n+1$ dimensional  Minkowski space  
with flat metric $\g = \text{diag}(1,1,\ldots -1)$. In what follows
$(\m,\h)$ will denote a Riemannian manifold with metric
$\h$; for simplicity we will restrict our attention to those manifolds
$(\m,\h)$ which are uniformly analytic; that is, the manifold can be covered
by a family of charts such that the metric and Christoffel symbol components
are analytic in each chart, with uniform exponential bounds on the
Taylor series coefficients.
Examples include $S^{\,m}$, $\R^{m}$, the hyperbolic plane, or any compact 
analytic manifold.

We are interested in maps
\begin{align}
\phi(x,t): \; (\rr^{n+1},\g) \longrightarrow (\m, \h). \label{mapping}
\end{align}
which are stationary with respect to compact variations of the Lagrangian
\begin{align}
\label{lagrangian}
L & = \int_{\rr^{n+1}} - \frac{1}{2} \text{Tr}_{\g} \phi^*(h) dv_{\g} \\
& =  \int_{\rr^{n+1}} - \frac{1}{2} \g^{\mu \nu} h_{\alpha \beta}(\phi)   
\frac{\partial \phi^{\alpha}}{\partial x^{\mu} }
 \frac{\partial \phi^{\beta}}{\partial x^{\nu} } dv_{\g} 
\label{lagrangian.local}
\end{align}
In \eqref{lagrangian.local} we have written \eqref{lagrangian} 
with respect to the coordinates $x^1, x^2, \ldots x^n, x^{n+1} = t$ on 
$\rr^{n+1}$ and local coordinates on $\m$. 
Stationary points of this Lagrangian are called wave maps,
and can be parametrized by the Cauchy problem
for the  wave map equation -  which is the 
Euler-Lagrange equation of \eqref{lagrangian}, 
\be{wavemap}
\bs
\Box \phi^k+
{\Gamma}_{\alpha \beta}^k(\phi) \partial_\mu\phi^\alpha\partial^{\mu}\phi^{\beta} 
& = 0.
\\
\phi[0] &= (f,g)
\es
\end{equation}
where $\Gamma^{k}_{\alpha \beta}$
are the Christoffel symbols
corresponding to the Riemannian metric $h$, and $\phi[T] = (\phi(T), \phi_t(T))$
denotes the Cauchy data of $\phi$ at time $T$.

A model to keep in mind is the case
with image $\m = S^{m-1} \subset \rr^m$, where the 
equations \eqref{wavemap} take the form (see e.g. \cite{struwe.barrett})
\be{wavemapsphere}
\bs
\Box \phi+ \phi \partial_\mu\phi^t \partial^{\mu}\phi
& = 0, \\
\phi[0] &= (f,g) \\
f(x) &\in S^{m-1}\\
f(x)^t g(x) &= 0
\es
\end{equation}
where we think of $\rr^m$ as an $m \times 1$ column vector, and
write $\phi^t$ for its transpose.
It is well
known that smooth solutions to \eqref{wavemap} will stay on the
sphere; the same result will hold for rough solutions by
a limiting argument assuming that the problem is well-posed in
the rough space.    

The equation \eqref{wavemap} is invariant under the scaling
\be{scaling-def}
\bs
\phi^\lambda(x,t) &= \phi(\frac x \lambda, \frac{t}\lambda)\\
f^\lambda(x) &= f(\frac x \lambda)\\
g^\lambda(x) &= \lambda^{-1} g(\frac x \lambda).
\es
\end{equation}
for any $\lambda > 0$. For data in $H^s \times H^{s-1}$, 
$s > n/2$, the initial value problem \eqref{wavemap} is called subcritical 
since $\|(f,g)\|_{\dot{H}^{s}\times \dot{H}^{s-1} }$ can be made small by
choosing $\lambda$ large in \eqref{scaling-def}.  A norm left
invariant by the scaling \eqref{scaling-def} is called critical.

%Global existence in all dimensions 
%from data in a small neighborhood of a point was shown in 
%\cite{kovalyov}, and from data in a small neighborhood of 
%a geodesic in \cite{sideris.harmonic}.  For large data in 
%dimensions $n \geq 3$, the problem is supercritical and
%smooth data has been shown to develop singularities. (\cite{shatah},
%\cite{shatah.shadi.blow}).  In the critical dimension $n = 2$, 
%global existence from large data is known if  assumptions on the 
%symmetries of the solution, and on the geometry of $\m$. (\cite{ 
 
In this paper we will usually restrict ourselves to the one-dimensional
case $n=1$,
where the analysis of the Cauchy problem 
\eqref{wavemap}
is simplified by introducing the null coordinates
\begin{align}
u & = x + t  \quad  v \; = \; x - t. \label{nullcoords} 
\end{align}
wherein the wave map equation is  
%\begin{align*}
%\partial_x & = \partial_u + \partial_v \quad   \partial_t \; = \;
%\partial_u - \partial_v \\
%\Box = -\partial_t^2 + \partial_x^2  & =  4 \partial_u \partial_v.
%\end{align*}
\be{wavemap.null}
\bs
\partial_u \partial_v \phi(u,v)  & = - \gamab (\phi) \cdot 
(\partial_u \phi^{\alpha} \partial_v \phi^{\beta})\\
\phi[0] &= (f,g)\\
\es
\end{equation}
and in the case of the sphere we have
\be{wave-map}
\bs
\phi_{uv} &= - \phi \phi_u^t \phi_v\\
\phi[0] &= (f,g)\\
f(x) &\in S^{m-1}\\
f(x)^t g(x) &= 0
\es
\end{equation}

We aim to show that for various sub-critical and critical
initial data spaces $\dd \subset C(\R)$, 
the Cauchy problem \eqref{wavemap.null} is 
locally and globally well posed in $\dd$, in the sense that
the solution operator exists and maps data from $\dd$ continuously into
a unique solution in  $C([0,T],\dd) \cap X$ for all $T > 0$, where
$X$ is some auxilliary space to be specified.
We also wish to show persistence of regularity, so that
a solution in a rough space $\tilde \dd$ whose initial data
is in a smooth space $\dd$ will stay in the smooth space $\dd$.
All the positive results are for the one-dimensional 
equation \eqref{wavemap.null}.
Finally, we complement these results with negative results.

For our local results in Theorem \ref{lwp-reg} below, 
we may assume that we are working 
in a single coordinate chart
since the norms used embed in the space of
H\"older continuous functions. (See section \ref{persist} below.)  

\begin{definition}\label{hsd-def} \cite{rauch.reed} Define the spaces
$H^{s,\delta}$ for $s, \delta \in \R$ by the norm
$$ \|\phi\|_{H^{s,\delta}} = \|\langle |\tau| + |\xi| \rangle^s \langle
|\tau| - |\xi| \rangle^\delta \hat\phi\|_{L^2_{\xi,\tau}}$$
where $\xi$, $\tau$ are the dual variables to $x$, $t$, and $\langle x\rangle
= (1 + |x|^2)^{1/2}$.
\end{definition}

The $H^{s,\delta}$  spaces first appeared in \cite{rauch.reed},
where Rauch-Reed study the propagation of singularities
of hyperbolic equations.  Our local 
well-posedness results are based on bilinear estimates in the 
$H^{s, \delta}$ spaces, versions of which first 
appeared in the higher dimensional 
work of \cite{bourg.I.II, kman.mach.nullform}.

\begin{theorem}\label{lwp-reg} (Local theory and persistence)
If\footnote{In the published version of this paper, local existence was claimed for $s>1/2$, but there was an error in that argument, pointed out to us by Kenji Nakanishi.  The full range of $s>1/2$ for local well-posedness has since been established in \cite{kenji}.} $s > 3/4$, then the Cauchy problem \eqref{wavemap.null} is
locally well-posed in $H^s \times H^{s-1}$ on some nontrivial time
interval $[0,T]$.  For any $\tilde{s} > 3/4$ 
one may choose $T$ to depend only 
on the $H^{\tilde s} \times H^{\stil -1}$ norm of the data.  
Furthermore
the solution to this problem is locally in $H^{s,s}$
with norm depending only on the $H^s \times H^{s-1}$ norm 
of the data.
\end{theorem}

In particular, if the $H^s$ solution cannot be continued past some
maximal time $T^*$, then the solution must blow up in 
$H^{\tilde s} \times H^{\tilde s-1}$ as $t \to T^*$
for all $\tilde s> 3/4$.  The last statement in Theorem \ref{lwp-reg}
will be needed for Theorem \ref{main} below.

\begin{theorem}\label{main}  (Global theory for Sobolev spaces) The
Cauchy problem \eqref{wavemap.null}
is globally well-posed for large data in $H^s$ for $s \geq 1$.  
In the case of the sphere, the Cauchy problem \eqref{wave-map}
is globally well-posed for large data in $H^s$ for
$1 > s > 3/4$.
\end{theorem}

In practice the distinction between functions which are globally in $H^s$ 
and those which are locally $H^s$ is unimportant, due to finite speed of
propagation.  It may be that Theorem \ref{main} can be extended to lower values of $s$ by a more sophisticated application of Lemma \ref{pohl_law} 
than is provided by our methods.  

\begin{theorem}\label{l11-main}  (Global theory for $L^{1,1}$)
The Cauchy problem \eqref{wavemap.null} is globally
well-posed with scattering for large data in the critical
space $L^{1,1}$ defined by
\begin{align} \| (f,g) \|_{L^{1,1}} & = \| f^\prime \|_{L^1} + \| f\|_{L^\infty} + \|g\|_{L^1}.
\label{l.11.def}
\end{align}
\end{theorem}

Our negative results are in Section \ref{negative}. 
We collect some previous observations on the ill-posedness of
superficially similar equations, and also show that the wave map problem
\eqref{wavemap} in certain coordinates is 
analytically ill-posed in the critical space $\dot H^{n/2}$ for $n \geq 1$.  
We also show there is no scattering when $n = 1$ 
without a suitable decay condition on the data.  Further ill-posedness
results in the one-dimensional case will appear in \cite{tao:ill}.

We now briefly discuss each of the positive results
and their relationship with previous literature. Careful 
surveys of regularity results and open
questions in spatial dimensions $n > 1$ can be found in
\cite{struwe.barrett},\cite{kman.barrett},\cite{grillakis.zurich},
and \cite{shatah.zurich}. 

The local theory of Theorem \ref{lwp-reg} is the $n=1$ version 
of higher-dimensional
results initiated in \cite{kman.mach.nullform}
and further studied in \cite{beals,grillakis.zurich, zhou.local, kman.mach.smoothing, 
kman.selberg, tataru,sogge}:  for $n \geq 2$, the wave map equation
is locally well-posed in the subcritical spaces 
$H^{s}(\rr^{n}), \, s > \frac{n}{2}$. Our methods are the same, 
but we have some additional simplifications due to the null-coordinates 
$u$, $v$ which are only available in one
dimension.  We prove the local results in Section \ref{persist}, 
after some abstract considerations in Section \ref{ab-persist}.

From Theorem \ref{lwp-reg} and energy conservation, one immediately obtains
Theorem \ref{main} for $s \geq 1$.  For $s < 1$ the energy conservation
law is not directly applicable, and to obtain our low regularity global
existence results we adapt some ideas of Bourgain \cite{bourg.2d} and 
a pointwise version of energy conservation observed
by Pohlmeyer\cite{pohl}.  This is the most 
involved part
of the paper, and occupies Sections \ref{prelim-sec}, 
\ref{wh-sec},
and \ref{proof-sec}.  Theorem \ref{main}  can almost certainly 
be extended to more general compact manifolds.  When $s \geq 2$, 
Theorem \ref{main}  
was previously shown in  \cite{gu,lady}, also \cite{ginebre,shatah}.  We have recently
learned that for $s \geq 1$, the result appears in \cite{zhou}.

%We note that the main interest of these global results is that
%there is no smallness condition on the data.  If the data is small in
%some $H^s$ norm, $s > \frac{1}{2}$, and hence small in the 
%Besov space $B^{\frac{1}{2},1}$, then global existence
%follows from the critical regularity results of 
%Tataru \cite{tataru, tataruII} after specializing
%the argument there to the one-dimensional setting.

Theorem \ref{l11-main} follows the usual pattern of global well-posedness
results for large-data in a critical space: one first proves
global existence for small data, then shows that the solution does
not concentrate.  Due to the simple structure and symmetries
of our equation and the data space, both of these steps are extremely 
elementary, especially when compared with other large-data critical results
e.g. \cite{grillakis.semilinear,bourg.critical,kap.energy,
struweshatah.energy, christ.spherical.wave, grillakis.equivariant,
shatah.zadeh.blow}.  Scattering is obtained by conformal
compactification.  
We prove this theorem in Section
\ref{l11}.

\noindent {\emph{ Acknowledgements:}}
Thanks to S. Klainerman and M. Machedon for explaining many ideas 
about the wave map equations, and to J. Bourgain for detailing his 
work in \cite{bourg.2d}.  Thanks to T. Sideris and D. Tataru
for timely advice, and to S.Y. Chang and L. Wang  for 
explaining ideas from \cite{changwangyang} which play a key 
role in Section \ref{global-h1-sec} below.  

\smallskip

\section{Abstract local existence and persistence theory}\label{ab-persist}

Questions of local existence
and persistence of regularity for nonlinear wave equations 
are very often handled by the method of Picard 
iteration, using estimates to control the nonlinearity.  For wave maps the
algebraic (or analytic) nature of the nonlinearity allows one to
formalize these estimates quite explicitly; this was done for instance
in Tataru \cite{tataru}.  In this section we describe the 
well-known abstract machinery which allows one to obtain existence 
and regularity from these estimates.
The statements of this section will be valid in every dimension $n \geq 1$.

We begin with the standard reduction of
local existence and persistence questions to estimates, which we set
abstractly.  We consider the Cauchy problem
\begin{equation}\label{abstract}
\begin{split}
L(\phi) &= N(\phi)\quad\quad t \in (-T_1,T_2)\\
\phi[0] &= f
\end{split}
\end{equation}
where $L$ is a linear evolution operator of order $d$, $N$ is a 
nonlinear operator such that $N(0)=0$, 
$0 \leq T_1, T_2 \leq \infty$ are times,
and the Cauchy data $f = (f_0, \ldots, f_{d-1})$ lies in some Banach space 
$\dd$.
We assume that a suitable coordinate chart has been chosen so that $\phi$
takes values in Euclidean space $\R^m$; this can be
done (locally at least) if $\dd$ embeds in the space of continuous functions.

We may rewrite this problem in integral form as
$$\phi = \eta (S(f) + L^{-1} N(\phi)),$$
where $S(f)$ is the solution to the homogeneous linear problem $L(\phi) = 0$
with Cauchy data $f$, $L^{-1} F$ is the solution to the
inhomogeneous problem $L \phi = F$ with Cauchy data $0$, and $\eta$
is any function which equals $1$ on $[-T_1, T_2]$.  For a rough
initial problem it will be advantageous to choose a smooth cutoff
$\eta$ (see \cite{bourg.I.II}).

We will always assume that the free problem is well-posed in $\dd$. 
For higher-dimensional wave equations this effectively 
restricts $\dd$ to the $L^2$-based family of spaces, but in one 
dimension many more spaces are available.

From the contraction mapping theorem we have the following local
existence metatheorem.  As this result is well-known, we omit
some details and rigor.  

Throughout the paper, we write $a \lesssim b$ to denote 
$a \leq C b$ for some large constant $C$.

\begin{lemma}[Local existence for small data] Let the notation be as above.
Suppose that there exists a reasonable\footnote{In fact, it's enough 
to assume that  $X$ can be densely approximated by test functions.}
Banach space $X$ of
functions in spacetime which obeys the estimates
\begin{align}
\Norm{\eta \sol(f)}_X &\lesssim \Norm{f}_{\dd} \label{data-X}\\
\Norm{\phi[T]}_{\dd} &\lesssim \Norm{\phi}_{X} \label{X-data}\\
\Norm{\eta L^{-1} (N(\phi) - N(\psi))}_{X} &\lesssim \
\Norm{\phi - \psi}_{X} (\Norm{\phi}_X + \Norm{\psi}_X)
\label{nl-est}
\end{align}
for all data $f$, $T \in (-T_1,T_2)$,
 and all spacetime functions $\phi, \psi$ with sufficiently small $X$ norm.

Then for sufficiently small $\epsilon$ depending only on
the constants in the above estimates, the Cauchy problem  
\eqref{abstract} is well posed in $(-T_1,T_2)$ for data $f \in \dd$,
with a unique solution in $X \cap C((-T_1,T_2),\dd)$, providing that 
$\|f\|_\dd \leq \epsilon$.  
\end{lemma}

\begin{proof}  If $\|f\|_D$ is sufficiently small, then the
assumptions imply the Picard iteration map
\begin{equation}\label{picard}
 \phi \mapsto \eta (S(f) + L^{-1} N(\phi))
\end{equation}
will be a contraction on a small neighborhood of the origin in $X$.
The contraction mapping theorem thus gives a unique solution on
this ball which depends continouously in $X$ on $S(f)$.
By \eqref{data-X} and \eqref{X-data} we get well-posedness.   Since
the solution is in $X$, it is in $L^\infty(\dd)$ by \eqref{X-data};
continuity in time follows from a straightforward approximation argument
using Schwartz functions.
\end{proof} 

One can relax the $(\Norm{\phi}_X + \Norm{\psi}_X)$ factor in the condition
\eqref{nl-est}, but we shall not need to do so in this paper.

A small modification of this argument allows one to get persistence of
regularity as long as the solution stays in a rougher space $\tilde D$
or $\tilde X$, providing of course that one has the appropriate
estimates.
\begin{lemma}[Persistence of regularity]
Assume $X \subset \xtil, \dd \subset \ddtil$ are spaces such that
$X$, $\dd$ satisfy \eqref{data-X} and \eqref{X-data}, and $\xtil$, $\ddtil$
satisfy \eqref{data-X}, and \eqref{nl-est}.  Assume
also that we have the estimate
\begin{equation}\label{x-xtil}
\begin{split}
\Norm{\eta L^{-1} (N(\phi) - N(\psi))}_{X} &\lesssim  \Norm{\phi - \psi}_{\xtil} 
\left(\| \phi \|_X + \| \psi\|_X\right)\\
& + \Norm{\phi - \psi}_{X}
\left( \Norm{\phi}_{\xtil} + \Norm{\psi}_{\xtil} \right)
\end{split}
\end{equation}
for all spacetime functions $\phi$, $\psi$ with sufficiently
small $\xtil$ norm.

Then there exists $\epsilon > 0$ such that
the problem \eqref{abstract} is well-posed in $(-T_1,T_2)$
for data $f \in \dd$, with unique solution in $X \cap C((-T_1,T_2),\dd)$,
providing that $\Norm{f}_{\ddtil} \leq \epsilon$.
\end{lemma} 

\begin{proof}
Define,
$$
\Norm{\phi}_Z = c(\eps \Norm{\phi}_{X} + C \Norm{\phi}_{\xtil})
$$
where $\eps = 1/\Norm{f}_{\dd}$ 
 and $c$, $C$ are large constants.
Then the assumptions imply that 
the map \eqref{picard} is a contraction on the unit ball in $Z$,
providing that $c$, $C$ are sufficiently large and $\epsilon$ is sufficiently 
small, hence the result.  
\end{proof}

One can show that the $X$-solution persists as long as the 
$\tilde X$ norm stays finite, but we shall not need that here.

We now specialize to the case of the wave map equation, in which
$L = \Box$ and 
$$N(\phi) = {\Gamma}_{\alpha \beta}^k(\phi) \partial_\mu\phi^\alpha\partial^{\mu}\phi^{\beta}.$$
We will assume that the data is small in $\dd$, and that $\dd$ is embedded
in the space $C(\R)$ of continuous data; this allows us to use a single chart
of coordinates.  In this case the required estimates can be simplified by the
identity
\begin{equation}\label{ident}
 \partial_\mu\phi \partial^{\mu}\psi = \frac{1}{2}[\Box(\phi \psi) - 
\phi \Box \psi - \psi \Box \phi].
\end{equation}
If $\Gamma$ depends polynomially on $\phi$, then $N(\phi) - N(\psi)$ 
can be decomposed by \eqref{ident}
 into quantities of the form $F \Box G$, where $F$
and $G$ are polynomials in $\phi$, $\psi$, and at least one of $F$, $G$
contains a factor of $\phi - \psi$.  Combining this with the previous
lemmas, one obtains

\begin{lemma}[Wave map local existence]\label{black-box}\cite{tataru} Suppose
that the Christoffel symbols $\Gamma(\phi)$ depend polynomially on $\phi$.
If $\dd \subset C(\R)$ and $X$ obey \eqref{data-X} and
\eqref{X-data}, as well as the estimates 
\begin{align}
\Norm{\phi \psi}_{X} &\lesssim   \Norm{\phi}_{X} 
\cdot \Norm{\psi}_{X}\label{algebra-property}\\
\Norm{\eta \Box^{-1} (\phi \Box \psi)}_{X} &\lesssim 
\Norm{\phi}_{X} \Norm{\psi}_{X}\label{box-property}
\end{align}
then the Cauchy problem \eqref{wavemap} is well-posed in $[-T_1,T_2]$
for data in $\dd$ provided that $\|f\|_\dd$ is sufficiently small.
\end{lemma}  

\begin{lemma}[Wave map persistence of regularity]\label{black.box.regularity}\cite{tataru} 
 Suppose
that the Christoffel symbols $\Gamma(\phi)$ depend polynomially on $\phi$.
Assume $X \subset \xtil, \dd \subset \ddtil \subset C(\R)$ are spaces such that
$X$, $\dd$ satisfy \eqref{data-X} and \eqref{X-data}, and $\xtil$, $\ddtil$
satisfy \eqref{data-X}, \eqref{algebra-property},
and \eqref{box-property}.  Assume
also that we have the estimates 
\begin{align}
\Norm{\phi \psi}_X &\lesssim \| \phi\|_X \|\psi\|_{\tilde X}
+ \| \psi\|_X \| \phi\|_{\xtil} \label{algebra-x-xtil}\\
\Norm{\eta \Box^{-1} (\phi \Box \psi)}_X &\lesssim 
\| \phi\|_X \|\psi\|_{\tilde X}
+ \| \psi\|_X \| \phi\|_{\xtil} \label{box-x-xtil}.
\end{align}
Then the Cauchy problem \eqref{wavemap} is well-posed in $[-T_1,T_2]$
for data in $\dd$ provided that $\|f\|_\ddtil$ is sufficiently small.
\end{lemma}

The same results hold if $\Gamma$ is uniformly analytic 
on the target manifold $\m$, since one can obtain the desired
estimates by expanding $\Gamma$ as a power series.  Note that the
geometry of $\m$ does not play any role in these results.

\section{Proof of Theorem \ref{lwp-reg}}\label{persist}

By finite speed
of propagation and the fact (from Sobolev embedding)
 that $H^s$ functions have some degree of
H\"older continuity for $s > 1/2$, we may assume 
that the data is compactly supported\footnote{For a proof that
$H^s$ functions can be localized, see Corollary \ref{stable}.}
 and stays within a single coordinate chart.  

Choose $\delta$ such that $3/4 < \delta < \stil, s, 1$.  
Since the $H^\stil \times
H^{\stil-1}$ norm of the data is bounded, the $H^\delta \times H^{\delta-1}$ 
norm is also bounded; we now show that by rescaling the data and shifting
coordinates we may make the $H^\delta \times H^{\delta-1}$ norm arbitrarily 
small.  (Our Sobolev norms are inhomogeneous
and do not obey an exact scaling identity, so one must take a little
care with this argument).

The equation \eqref{wavemap.null} is invariant under the scaling
\eqref{scaling-def}. Eventually we will choose $\lambda$ depending
only on $\delta$, $\stil$, and the $H^\stil \times H^{\stil-1}$ norm
of the data.
Thus to obtain well-posedness for the original data up to time $1/\lambda$
it suffices to get well-posedness up to time $1$ for the data $f^\lambda$,
$g^\lambda$.  By finite speed of propagation we may restrict $f^\lambda$, 
$g^\lambda$ to an interval of length $4$ centered at some $x_0$; 
by translation invariance we may make $x_0 = 0$. 
By shifting the origin of the coordinate system we may replace $f^\lambda$ by 
$f^\lambda - c$ for some constant $c$.  In particular, it 
suffices to get 
well-posedness up to time $1$ for the data
\begin{align*}
\tilde f(x) &= \chi(x) (f(\frac{x}{\lambda}) - \overline{f})\\
\tilde g(x) &= \lambda^{-1} \chi(x) g(\frac{x}{\lambda}) 
\end{align*}
where $\chi$ is a standard compactly supported cutoff function,
and $\overline{f} = \int f(\frac{x}{\lambda}) \psi(x)\ dx$,
where $\psi$ is a standard bump function with unit mass supported
near $\chi$.

We claim that we may make the $H^\delta \times H^{\delta-1}$ norm of
$(\tilde f, \tilde g)$ arbitrarily small, by choosing $\lambda$
sufficiently large (but depending only on the $H^\stil \times H^{\stil-1}$
norm of the original data).  More precisely,

\begin{lemma}\label{rellich}  If $\delta < \stil$ and $\delta \leq 1$, we have
\be{rescale}
\| \tilde f \|_{H^{\delta}} + \| \tilde g \|_{H^{\delta - 1}} 
\lesssim \lambda^{-\eps} ( \|f\|_{H^\stil} + \|g\|_{H^{\stil - 1}} )
\end{equation}
for all $\lambda \gg 1$, where $\eps > 0$ is a small number depending 
only on $\delta$, $\stil$.
\end{lemma}

\begin{proof}
The contribution of $g$ is easily handled by the rescaling properties
of $H^{\delta-1}$:
$$ \| \tilde g \|_{H^{\delta - 1}} \lesssim \lambda^{-1} \|  g(\frac x \lambda) \|_{H^{\delta -1}} \lesssim \lambda^{-1} \lambda^{\frac{1}{2}} 
\| g \|_{H^{\delta-1}} \lesssim \lambda^{-1/2} \| g \|_{H^{\stil - 1}},$$
so we may assume that $g=0$.  It suffices
to check the cases when the Fourier transform of $f$ are supported on
$| \xi | \lesssim \lambda$ and $| \xi | \gtrsim \lambda$.

We first consider the case when $|\xi| \lesssim \lambda$.  Since $\delta \leq 1$
and $\tilde f$ is compactly supported, the $H^\delta$ norm is controlled by
the $C^1$ norm, and so it suffices to control the quantity
$$ \| \chi (f(\frac{x}{\lambda}) - \overline{f}) \|_{C^1}.$$
But a computation shows that this is majorized by $\lambda^{-1} \|f\|_{C^1}$,
which by Sobolev embedding is majorized 
by $\lambda^{-1} \|f\|_{H^{3/2 + \eps}}$, which by the frequency support
assumption on $f$ is controlled by $\lambda^{1/2 - \stil + \eps}
\|f\|_{H^\stil}$, which gives the desired estimate if $\eps$ is chosen
sufficiently small.

We now consider the case when $|\xi| \gtrsim \lambda$.  In this case
we use the triangle inequality to estimate
$$ \| \tilde f\|_{H^\delta} \lesssim \|\chi f(\frac{x}{\lambda}) \|_{H^\delta}
+ \| \chi \overline{f} \|_{H^\delta} \lesssim \| f(\frac{x}{\lambda}) 
\|_{H^\delta} + |\overline{f}|.$$
By the frequency support assumption on $f$ and Plancherel's theorem, we have
$$ \| f(\frac{x}{\lambda}) \|_{H^\delta} \sim
\| f(\frac{x}{\lambda}) \|_{\dot H^\delta} \sim \lambda^{1/2 - \delta}
\| f\|_{\dot H^\delta} \lesssim \lambda^{1/2 - \delta} \|f\|_{H^\stil}$$
as desired.  To control $\overline{f}$, we use Plancherel's theorem to write
$$ \overline{f} = \int \lambda \hat{f}(\lambda \xi) \hat{\psi}(\xi)\ d\xi
= \int \hat{f}(\xi) \hat{\psi}(\frac{\xi}{\lambda})\ d\xi.$$
From the support hypothesis on $\hat{f}$ and Cauchy-Schwarz, this is
estimated by
$$ (\int_{|\xi| \gtrsim \lambda} |\hat{f}(\xi)|^2 \langle \xi\rangle^{2\stil}\ d\xi)^{1/2} (\int_{|\xi| \gtrsim \lambda} \hat{\psi}(\frac{\xi}{\lambda}) \langle \xi\rangle^{-2\stil}\ d\xi)^{1/2};$$
since $\psi$ is rapidly decreasing, this is majorized by $\lambda^{1/2 - \stil}
\|f\|_{H^{\stil}}$, which is acceptable.
\end{proof}

It is likely that a version of the above lemma can also be proven by
Rellich's lemma and a compactness argument using
the nonconcentration of $H^\delta$ norm for smooth functions, but we 
shall not do so here.

To finish
the proof of Theorem \ref{lwp-reg}, we have to show that the equation
\eqref{wavemap.null} is locally well-posed in $H^s$ up to time 1 with a 
solution in $H^{s,s}$, whenever the $H^\delta \times H^{\delta-1}$ norm
of the data is sufficiently small.  

We apply Lemma \ref{black.box.regularity}, with $\dd = H^s \times H^{s-1}$,
$\tilde \dd = H^{\delta} \times H^{\delta - 1}$, $X = H^{s,\delta}$,
$\tilde X = H^{\delta,\delta}$, $T_1=T_2=1$, and a smooth cutoff
$\eta$.  Assuming we can verify all the estimates in the lemma,
this gives well-posedness in $H^s$ up to time $1$ with a solution
in $H^{s,\delta}$.  At the end of this section we shall improve
this to $H^{s,s}$.

Of course, it still remains to verify the hypotheses in Lemma 
\ref{black.box.regularity}.  More precisely, 
we need to show that $X$, $D$ satisfy \eqref{data-X},
\eqref{X-data}, that $\xtil$, $\ddtil$ satisfy \eqref{data-X}, 
\eqref{algebra-property}, and \eqref{box-property}, and that
\eqref{algebra-x-xtil}, \eqref{box-x-xtil} hold; the
inclusion $\dd \subset C(\R)$ follows from Sobolev embedding.

We first take advantage of the null coordinates to
rewrite the $H^{s,\delta}$ norms in terms of
product Sobolev spaces $H^{s_1}_u H^{s_2}_v = H^{s_2}_v H^{s_1}_u$ defined by
$$ \| \phi \|_{H^{s_1}_u H^{s_2}_v} = \| D_u^{s_1} D_v^{s_2} \phi \|_{L^2_{u,v}},$$
where $D_u$ and $D_v$ are the Fourier multipliers corresponding to
$\langle \mu\rangle$, $\langle \nu \rangle$ respectively, and $\mu, \nu$
are the frequency variables dual to $u,v$.  We define the 
one-dimensional Sobolev spaces $H^s_u$, $H^s_v$ in the usual manner.

By Plancherel's theorem one can easily verify that
\begin{equation}\label{product}
H^{s,\delta} = H^s_u H^\delta_v \cap H^s_v H^\delta_u
\end{equation}
when $\delta \leq s$.
Thus to prove estimates concerning the $H^{s,\delta}$ spaces in
$\R^{1+1}$, it suffices to prove estimates on product Sobolev spaces.
We collect the estimates we will need below, and then use them 
to prove the requirements of Lemma \ref{black.box.regularity}.

We first begin with a standard result regarding multiplication of
one-dimensional Sobolev spaces; we will use variants of this argument
in other places in this paper.

\begin{lemma}\label{easy-mult} If $s,s^\prime$ are real numbers such that
 $s > 1/2$ and $s \geq s^\prime \geq -s$, then
for all test functions $\phi$, $\psi$
$$ \| \phi \psi \|_{H^{s^\prime}_u}
\lesssim \| \phi \|_{H^{s}_u} \| \psi \|_{H^{s^\prime}_u}.$$
\end{lemma}

\begin{proof}
We may assume that the norms on the right-hand side are equal to one.
By Plancherel's theorem it suffices to show that
\begin{equation}\label{ft}
 \langle \mu \rangle^{s^\prime} (\hat \phi * \hat \psi)(\mu)
= \int_{\mu_1 + \mu_2 = \mu} \langle \mu_1 + \mu_2 \rangle^{s^\prime}
\hat \phi(\mu_1) \hat \psi(\mu_2)\ d\mu_1
\end{equation}
is in $L^2_\mu$.
Since the right-hand side norms depend only on the size
of $\hat \phi$, $\hat \psi$, we may assume that these functions
are non-negative.

We observe the elementary inequality
$$
\langle \mu_1 + \mu_2 \rangle^{s^\prime}
\lesssim \langle \mu_2 \rangle^{s^\prime} + \langle \mu_1 \rangle^{s}
\langle \mu_2 \rangle^{s^\prime - s} + \langle \mu_1 + \mu_2 \rangle^{-s}
\langle \mu_1 \rangle^s \langle \mu_2 \rangle^{s^\prime},
$$
which is easily shown by checking the cases $\langle \mu_1 \rangle
\ll \langle \mu_2 \rangle$, $\langle \mu_1 \rangle \gg  \langle \mu_2 \rangle$,
$\langle \mu_1 \rangle \sim  \langle \mu_2 \rangle$ seperately.  By
applying this estimate to \eqref{ft} 
and using Plancherel's theorem
again, we see that it suffices to show that
$$\phi (D_u^{s^\prime} \psi),  (D^s_u \phi) (D^{s^\prime - s}_u \psi),
D^{-s}_u [ (D^s_u \phi) (D^{s^\prime}_u \psi) ]$$
are each in $L^2_u$.

The first function is a product of an $H^s_u$ and an $L^2_u$ function,
and is thus in $L^2$ by the Sobolev embedding $H^s_u \subset L^\infty_u$.
The second function is a product of an $L^2_u$ and an $H^s_u$ function
and is treated similarly.  To show that the last function is in $L^2$,
it suffices by the Sobolev embedding $D^{-s}_u L^1_u \subset L^2_u$
to show that $(D^s_u \phi) (D^{s^\prime}_u \psi)$ is in $L^1_u$.
But this follows from H\"older's inequality since the two factors
are in $L^2_u$.  Thus $D_u^{s^\prime} (\phi \psi)$ is in $L^2_u$ as
desired.
\end{proof}

The same argument applies of course to the $v$ variable. Working
in both the $u$ and $v$ variables we obtain, 
%\footnote{One has to take a little care in deciding the order in which to place mixed
%Lebesgue norms.}, 

\begin{lemma}\label{product-h1}  If $s_1, s_2 > 1/2$, and
$s_1 \geq s_1^\prime \geq -s_1$, $s_2 \geq s_2^\prime \geq -s_2$,
then for all test functions $\phi$, $\psi$
\begin{align}
\| \phi \psi \|_{H^{s_1^\prime}_u H^{s_2^\prime}_v}
& \lesssim \| \phi \|_{H^{s_1}_u H^{s_2}_v} 
\| \psi \|_{H^{s_1^\prime}_u H^{s_2^\prime}_v}. \label{ph1-a} \\
\| \phi \psi \|_{H^{s_1^\prime}_u H^{s_2^\prime}_v}
& \lesssim \| \phi \|_{H^{s_1^\prime}_u H^{s_2}_v} 
\| \psi \|_{H^{s_1}_u H^{s_2^\prime}_v}. \label{ph1-b} 
\end{align}
\end{lemma}

\begin{corollary}\label{stable} The one-dimensional and product Sobolev spaces
are stable under multiplication by bump functions.  In particular,
if $\eta$ is a bump function and $\tilde \eta$ is a Schwarz function which
is non-zero on the support of $\eta$, then
$\eta \phi \in H^{s_1}_u H^{s_2}_v$ whenever
$\tilde \eta \phi \in H^{s_1}_u H^{s_2}_v$.
\end{corollary}

Finally, we need the following lemma on the 
smoothing properties of $\Box^{-1}$.  For previous instances of 
this lemma in higher dimensions and for differential operators 
other than $\Box$, see \cite{bourg.I.II,kpv.negative,kman.mach.smoothing}.

\begin{lemma}\label{integ}  If $\eta$ is a fixed bump function and
$s_1, s_2 \geq 1/2$ with $s_1+s_2 > 3/2$ and $|s_1-s_2| \leq 1$, then
$$ \| \eta \Box^{-1} \phi \|_{H^{s_1}_u H^{s_2}_v} \leq C_\eta
\| \phi \|_{H^{s_1-1}_u H^{s_2-1}_v}$$
for all test functions $\phi$.
\end{lemma}

\begin{proof} See \cite[Lemma 2.5]{kenji}.  (An argument in the published version of this paper omitted the necessary conditions $s_1+s_2>3/2$ and $|s_1-s_2| \leq 1$, and were incorrect; this is the reason why our local well-posedness results are restricted to $s>3/4$ rather than $s>1/2$.  We thank Kenji Nakanishi for pointing out the issue, which is further discussed in \cite{kenji}.)
\end{proof}

We can now prove the estimates necessary to apply 
Lemma  \ref{black.box.regularity}.

We first prove \eqref{data-X}, which in this context is
\be{dx-wave}
\| \eta S(f,g) \|_{H^{s,\delta}} \lesssim \|f\|_{H^s} + \|g\|_{H^{s-1}}.
\end{equation}
We observe that $S(f,g)$ can be written as 
$F(u) + G(v)$ for some compactly supported
$H^s$ functions $F$, $G$.  By Corollary
\ref{stable}
it thus suffices to show that 
$F(u) \eta(v)$ and $G(v) \eta(u)$ are in
$H^{s,\delta}$ for one-dimensional cutoff functions $\eta$.  But this 
follows from \eqref{product}.  A similar argument shows that
$\xtil$, $\ddtil$ also obey \eqref{data-X}.

We next prove \eqref{X-data}, which in this context is
\begin{equation}\label{xd-wave}
 \| \phi(T) \|_{H^s} + \| \phi_t(T) \|_{H^{s-1}}
\lesssim \| \phi \|_{H^{s,\delta}}.
\end{equation}
It suffices to show that $(D^s \phi)(T)$ is in $L^2$ for all
multipliers $D^s$ which are symbols of order $s$.  Since
the symbol of $D^s$ is majorized by that of $D_u^s + D_v^s$,
we can thus decompose $D^s \phi = \psi_1 + \psi_2$, where
$\psi_1$ and $\psi_2$ are in $L^2_u H^\delta_v$ and $L^2_v H^\delta_u$ 
respectively (by \eqref{product}).
The claim then follows from Sobolev embedding and the fact that
the $t=T$ trace of an $L^2_u L^\infty_v$ or $L^2_v L^\infty_u$
function is in $L^2$.

We now prove \eqref{algebra-x-xtil} and \eqref{box-x-xtil};
the proof that $\xtil$, $\ddtil$ satisfy
\eqref{algebra-property} and \eqref{box-property} will
follow by specializing the following arguments (which do not
need the hypotheses $\delta < 1,s$) to the case $s = \delta$.
In our context, the estimates to prove are
\begin{align}
 \| \phi \psi \|_{H^{s,\delta}} &\lesssim \| \phi \|_{H^{s,\delta}} 
\| \psi\|_{H^{\delta,\delta}} + 
\| \phi \|_{H^{\delta,\delta}} \| \psi\|_{H^{s,\delta}}\label{algebra}\\
 \| \eta \Box^{-1} (\phi \Box \psi) \|_{H^{s,\delta}} &\lesssim \| \phi \|_{H^{s,\delta}}
\| \psi\|_{H^{\delta,\delta}}
+ \| \phi \|_{H^{\delta,\delta}} \| \psi\|_{H^{s,\delta}}.\label{box}
\end{align}

The estimate \eqref{algebra} follows immediately from
\eqref{product} and \eqref{ph1-a}, so it only remains to show
\eqref{box}.
It suffices by \eqref{product}
and $u$-$v$ symmetry to estimate the $H^s_u H^{\delta}_v$ norm of
$\eta \Box^{-1} (\phi \Box \psi)$, which by
Lemma \ref{integ} is controlled by the $H^{s-1}_u H^{\delta - 1}_v$ norm of
$\phi D_u D_v \psi$.  We now divide into two cases.  If $s-1 \leq \delta$,
then by \eqref{ph1-a} (with $\psi$ replaced by $D_u D_v \psi$) we have
$$ \| \phi D_u D_v \psi \|_{H^{s-1}_u H^{\delta - 1}_v}
\lesssim \| \phi \|_{H^{\delta}_u H^{\delta}_v} 
\| D_u D_v \psi \|_{H^{s-1}_u H^{\delta - 1}_v},$$
which gives \eqref{box}.  When $s-1 > \delta$, the proof is similar
but \eqref{ph1-a} is replaced by the following lemma (with $\psi$
replaced by $D_u D_v \psi$):

\begin{lemma} If $s-1 > \delta > 1/2$, then
$$ \| \phi \psi \|_{H^{s-1}_u H^{\delta - 1}_v}
\lesssim \|\phi\|_{H^{s}_u H^{\delta}_v} \|\psi\|_{H^{\delta-1}_u H^{\delta-1}_v}
+ \|\phi\|_{H^{\delta}_u H^{\delta}_v} \|\psi\|_{H^{s-1}_u H^{\delta-1}_v}.
$$
\end{lemma}

\begin{proof}  We repeat the argument in Lemma \ref{easy-mult}.  It suffices
to estimate
\be{exp}
\int_{\mu_1 + \mu_2 = \mu} \int_{\nu_1 + \nu_2 = \nu}
\langle \mu_1 + \mu_2 \rangle^{s-1} \langle \nu_1 + \nu_2 \rangle^{\delta-1}
\hat \phi(\mu_1,\nu_1) \hat \psi(\mu_2,\nu_2)\ d\mu_1 d\mu_2
\end{equation}
in $L^2_\mu L^2_\nu$,
and we may assume as before that $\hat \phi$, $\hat \psi$ are non-negative.

By Plancherel's theorem and the easily verified inequalities
\bas
\langle \mu_1 + \mu_2 \rangle^{s-1}
&\lesssim \langle \mu_1 \rangle^s \langle \mu_2 \rangle^{-1}
+ \langle \mu_2 \rangle^{s-1}\\
\langle \nu_1 + \nu_2 \rangle^{\delta - 1}
&\lesssim \langle \nu_2 \rangle^{\delta - 1} + \langle \nu_1 \rangle^{\delta}
\langle \nu_2 \rangle^{- 1} + \langle \nu_1 + \nu_2 \rangle^{-\delta}
\langle \nu_1 \rangle^\delta \langle \nu_2 \rangle^{\delta-1},
\end{align*}
the $L^2_\mu L^2_\nu$ norm of \eqref{exp} is majorized by the $L^2_u L^2_v$
norms of
\bas
(D^s_u \phi) (D^{-1}_u D^{\delta-1}_v \psi),\quad
(D^s_u D_v^\delta \phi) (D^{-1}_u D^{-1}_v \psi),&\quad
D_v^{-\delta} [(D^s_u D^\delta_v \phi) (D^{-1}_u D^{\delta-1}_v \psi)],\\
(\phi) (D^{s-1}_u D^{\delta-1}_v \psi),\quad
(D_v^\delta \phi) (D^{s-1}_u D^{-1}_v \psi),&\quad
D_v^{-\delta} [(D^\delta_v \phi) (D^{s-1}_u D^{\delta-1}_v \psi)].
\end{align*}
The $L^2_u L^2_v$ norms of the
first three expressions are controlled by the $H^s_u H^\delta_v$
norm of $\phi$ and the $H^{\delta-1}_u H^{\delta-1}_v$ norm of $\psi$,
using the H\"older and Sobolev inequalities (in particular, the
fact that the product of an $L^2_u H^\delta_v$ and an $L^2_v H^\delta_u$
function is in $L^2_u L^2_v$) as in the proof of Lemma 
\ref{easy-mult}.  The last three expressions are similarly controlled
by the $H^{\delta}_u H^\delta_v$ norm of $\phi$ and the $H^{s-1}_u
H^{\delta-1}_v$ norm of $\psi$.
\end{proof}

Finally, we show that the solution $\phi$ is locally in $H^{s,s}$.  From the
above we have that $\phi$ is locally in $H^s_u H^\delta_v$.  Since
this space is an algebra by Lemma \ref{easy-mult} and $\Gamma$
is analytic, we see that $\Gamma(\phi) \in H^s_u H^\delta_v$.
Also we have $\phi_u \in H^{s-1}_u H^\delta_v$, while
a symmetrical argument gives $\phi_v \in H^\delta_u H^{s-1}_v$.
We now divide into 
 the cases $\delta \geq s-1$ and $\delta < s-1$.  If
$\delta \geq s-1$ then \eqref{wavemap.null} and Lemma \ref{easy-mult} now gives
$$ \phi_{uv} = \Gamma(\phi) \phi_u \phi_v \in 
(H^s_u H^\delta_v)(H^{s-1}_u H^\delta_v)(H^\delta_u H^{s-1}_v) \subset 
H^{s-1}_u H^{s-1}_v \hbox{ (locally)}$$
Since $\eta \phi = \eta S(f,g) + \eta \Box^{-1} \phi_{uv}$, the claim
then follows from \eqref{dx-wave} (with $\delta=s$) and Lemma \ref{integ}.

If $\delta < s-1$, then the above argument will only yield that
$\phi_{uv}$ is in $H^{s-1}_u H^\delta_v \cap H^\delta_u H^{s-1}_v$,
so that $\phi$ is in $H^{s,\delta+1}$.  One then iterates
the above argument, with $\delta$ replaced by $\delta+1$, until one
eventually obtains $H^{s,s}$ control on $\phi$.

\section{A pointwise conservation law, and consequences}\label{global-h1-sec}

In this section we introduce a pointwise conservation law for the 
one-dimensional wave map equation which is 
special to the one-dimensional case.  This law was first 
observed by Polhmeyer\cite{pohl}.  (See also \cite{shatah}.) An 
identity key to our work here, \eqref{hardwire}, is motivated by 
\cite{changwangyang}.  

\begin{lemma}\label{pohl_law}  If $\phi$ is a smooth solution to
\eqref{wavemap}, then
the quantity $|\phi_u|_\h$ is constant with respect to $v$, and the quantity $|\phi_v|_\h$ is
constant with respect to $u$, where we use $|x|_\h$ to denote
the length of a tangent vector $x$ in $\m$ with respect to the Riemannian
metric $\h$.
\end{lemma}

\begin{proof}
The energy-momentum tensor $T_{\alpha \beta}$ for wave maps is 
\begin{align}
\label{energy-momentum}
T_{\alpha \beta} & = \frac{1}{2} \big( \left< \partial_{\alpha} \phi,
\partial_{\beta} \phi \right>  - \frac{1}{2} \g_{\alpha \beta}
\left< \partial_{\mu} \phi, \partial^{\mu} \phi \right> \big)
\end{align}
where $\left<\phi, \psi \right> \; = \h_{\mu \nu} \phi^{\mu} \psi^{\nu}$
is the inner product on $\m$.  Recall that in all dimensions, the tensor $T$ is 
divergence free,
\begin{align}
\partial^{\alpha} T_{\alpha \beta} & = 0.
\label{divfree}
\end{align}
In $\rr^{1+1}$, $\g_\alpha^{\;\alpha} = 2$ and so $T$ is also trace free,
$$ T_{\alpha}^{\; \alpha} = \frac{1}{2} \big( \left< \partial_{\alpha} \phi,
\partial^{\alpha} \phi \right> -  \frac{1}{2} \g_\alpha^{\; \alpha}
\left< \partial_{\mu} \phi, \partial^{\mu} \phi \right> \big) = 0.$$

We write these two facts in null coordinates $u$, $v$.  The 
trace-free property gives $T_{uv} + T_{vu} = 0$; since $T$ 
is symmetric we thus have $T_{uv} = T_{vu} = 0$.
The divergence-free property then gives
$$ \partial_v T_{uu} = \partial_u T_{vv} = 0,$$
so that $T_{uu}$ is constant with respect to $v$, and $T_{vv}$ is constant with respect
to $u$.  The claim then follows since $T_{uu} = \frac{1}{2}|\phi_u|_\h^2$,
and $T_{vv} = \frac{1}{2}|\phi_v|_\h^2$.
\end{proof}

Although this lemma is phrased for smooth solutions, the result 
extends to rough solutions by applying a limiting argument and using 
the local well-posedness theory from Theroem \ref{lwp-reg}.  
Note that Lemma \ref{pohl_law} obviously holds as well for solutions of 
the free wave equation.

% (this is equivalent to trace-free ???) One can also obtain these 
%identities from Noether's theorem and the invariance
%of the wave map equation under the 
%conformal transformation $u \to \Phi(u)$, $v \to \Psi(v)$
%for arbitrary diffeomorphisms $\Phi$, $\Psi$.

In the case when the target manifold is a sphere, a more direct proof
is available.  Since the solution $\phi$ is on the sphere, we
have $\phi^t \phi = 1$.  Differentiating with respect to
$u$ we obtain $\phi^t \phi_u = 0$.  Combining this with \eqref{wave-map}
we obtain the useful identity
\be{hardwire}
\phi_{uv} = R \phi_u,
\end{equation}
where $R$ is the anti-symmetric matrix
\be{r-def}
R = \phi_v \phi^t - \phi \phi_v^t.
\end{equation}
The anti-symmetry of $R$ implies that $|\phi_u|^2$ is constant in the $v$ 
direction:
\be{asym}
 \partial_v |\phi_u|^2 = 2 \phi_u^t \phi_{uv} = 2 \phi_u^t R \phi_u = 0,
\end{equation}
and the other conservation law in Lemma \ref{pohl_law} is proven similarly.

Lemma \ref{pohl_law} can be viewed as a pointwise form of energy conservation,
and has many consequences.  For $H^1$ solutions it implies the estimates
\begin{align}
\| \phi_u \|_{L^2_u L^\infty_v} &\lesssim \|f\|_{\dot H^1} + \|g\|_2 \label{h1-u}\\
\| \phi_v \|_{L^2_v L^\infty_u} &\lesssim \|f\|_{\dot H^1} + \|g\|_2, \label{h1-v}
\end{align}
which in turn show that the $H^1 \times L^2$ norm\footnote{Of course,
one can show energy conservation much more directly, but the above approach
is more robust, and can be extended to regularities below the $H^1$ norm.}
 of the solution $\phi[t]$
at time $t$ is bounded uniformly in $t$.  Combining this with
Theorem \ref{lwp-reg} one obtains Theorem \ref{main}
for $s \geq 1$.

From \eqref{wavemap}, \eqref{h1-u}, \eqref{h1-v}, and the 
assumption that $\Gamma(\phi)$ is uniformly
bounded, we get the spacetime estimate
\be{h1-uv}
\| \phi_{uv} \|_{L^2_u L^2_v} \lesssim (\|f\|_{\dot H^1} + \|g\|_2)^2. 
\end{equation}

The following Corollary to Lemma \ref{pohl_law} states
that when the initial data is essentially compactly supported, 
the solution to \eqref{wave-map} resolves to an exact free solution
in finite time.
\begin{corollary}\label{scatter}  Suppose that $H^1$ Cauchy data $(F, G)$ 
are given such that $F^\prime$, $G$ are supported on the interval
$[-T,T]$.  Then the global solution $\Phi(u,v)$ to \eqref{wave-map} with
this data is
constant on the quadrants $[T,\infty) \times [T,\infty)$, $[T,\infty)
\times (-\infty,-T]$, $(-\infty,-T] \times [T,\infty)$, 
and $(-\infty,-T] \times (-\infty,-T]$, is constant in the $v$
direction on the strips $[-T,T] \times [T,\infty)$, $[-T,T] \times (-\infty,T]$,
and is constant in the $u$ direction in the strips $[T,\infty) \times [-T,T]$,
$(-\infty,T] \times [T,T]$.
\end{corollary}

\begin{figure}[htbp] \centering
\ \psfig{figure=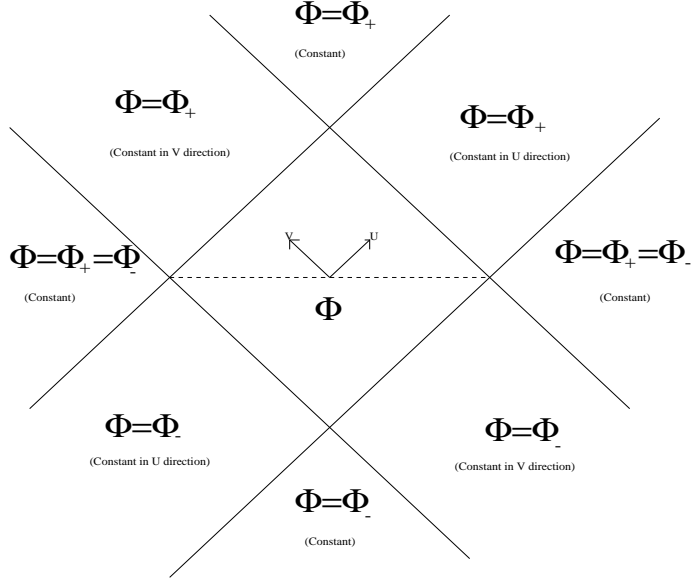,height=3in,width=3.6in}
\caption{Scattering of solutions $\Phi$ with essentially compactly supported
data.}
        \label{fig:diamond}
        \end{figure}

In particular, we see that $\Phi$ scatters exactly
to a free solution $\Phi_+$ when $t > T$ and to another free solution $\Phi_-$
when $t < -T$.  (See Figure \ref{fig:diamond}).

\section{Global existence in $H^s$, $3/4 < s < 1$: preliminaries}\label{prelim-sec}

We now turn to the second claim in Theorem \ref{main}.  Fix
$3/4 < s < 1$.  We have to show that the Cauchy problem
\eqref{wave-map} for the sphere is globally well-posed
for data which is locally in $H^s$.  It suffices to show local
well-posedness on some time interval $(-T_0,T_0)$, where we fix $T_0$ to be
an arbitrary large time.  By finite speed of propagation we may
assume that the data becomes constant outside of the interval $[-4T_0,4T_0]$.
In particular, we have
\be{compact}
f^\prime, g \hbox{ is supported on } [-4T_0,4T_0].
\end{equation}
We also make the a priori assumption that the data is in $H^1$;
this assumption will be removed by a density argument.  More 
precisely, we will assume that
\begin{align}
 \| f\|_{H^s} + \|g\|_{H^{s-1}} &\leq C_0\label{hs-est}\\
 \| f\|_{H^1} + \|g\|_{L^2} &\leq M \label{h1-finite}
\end{align}
where $C_0, M$ are arbitrary constants.  Henceforth all constants
will be allowed to depend on $C_0$, \emph{but not on} $M$.  We will
use the quantities $C$, $N$, $\eps$ to denote positive
constants that vary from line to line.

Since the data is in $H^1$, there is a unique global $H^1$ solution
$\phi$.  We aim to show the $H^s$ norm of the solution 
is bounded by a quantity which depends polynomially on $T_0$ but
is independent of $M$:
\be{bs-cons}
\| \phi(t) \|_{H^{s}} + \| \phi_t(t)\|_{H^{s-1}} \lesssim 
T_0^N \hbox{ for all } |t| \leq T_0.
\end{equation}
Then by Theorem \ref{lwp-reg} and a limiting argument, the same
estimate holds without the condition \eqref{h1-finite}, and
one obtains well-posedness in $H^s$ on the interval $(-T_0,T_0)$.

It remains to prove \eqref{bs-cons}.  When $s=1$ such
an estimate can be obtained from \eqref{h1-u} and \eqref{h1-v},
so it is natural to look for variants of \eqref{h1-u},
\eqref{h1-v} (and perhaps \eqref{h1-uv}) which apply for data which are 
rougher than
$H^1$.  

For the free equation $\Box \phi = 0$, conservation
of the $H^s$ norm for $s<1$ is shown by applying fractional integration
operators to the equation and then applying the energy conservation law.  
Thus a first guess might be to apply the operators $D_u^{s-1}$, $D_v^{s-1}$
to the above equations.  The fact that these operators
(for fractional $s$) are not local is inconvenient for technical reasons,
and we will instead apply the following modified fractional integration 
operators.
\begin{definition}\label{dtil}  For any $s \in \R$, let $m_s(\xi)$ be
the convolution of $\langle \xi\rangle^s$ with $\hat \eta$, where
$\eta(x)$ is a bump function with non-negative Fourier transform.  We let
$\tilde D_u^s$, $\tilde D_v^s$ be the Fourier multipliers
corresponding to $m_s(\mu)$, $m_s(\nu)$ respectively.
\end{definition}

These operators behave like the usual fractional differentiation and
integration operators, but have a compactly supported kernel.
Note that $m_s(\xi)$ is comparable to $\langle \xi\rangle^s$, so
one can replace $D_u^s$, $D_v^s$ by $\tilde D_u^s$, $\tilde D_v^s$
respectively in the definitions of the Sobolev spaces defined earlier.
That the $\tilde D_u^s$ operators are not perfectly multiplicative in $s$ is
irrelevant for our purposes.

If one informally pretends that $\tilde D_u^{s-1}$, $\tilde D_v^{s-1}$
commute with the wave map equation, then \eqref{h1-u}, \eqref{h1-v}, 
\eqref{h1-uv} informally yield
\begin{align}
\| \phi_u \|_{H^{s-1}_u L^\infty_v} &\lesssim \|f\|_{H^s} + \|g\|_{H^{s-1}} \label{hs-u}\\
\| \phi_v \|_{H^{s-1}_v L^\infty_u} &\lesssim \|f\|_{H^s} + \|g\|_{H^{s-1}} 
\label{hs-v}\\
\| \phi_{uv} \|_{H^{s-1}_u H^{s-1}_v} &\lesssim (\|f\|_{H^s} + \|g\|_{H^{s-1}})^2, \label{hs-uv}
\end{align}
where the spaces $H^{s-1}_u L^\infty_v$, $H^{s-1}_v L^\infty_u$
are defined\footnote{For technical reasons caused by the $L^\infty$ norm,
one has to take some care in defining these spaces; for instance, one cannot
simply replace $\tilde D_u^{s-1}$ by $D_u^{s-1}$.} 
\be{hl-def}
\| \phi \|_{H^{s-1}_u L^\infty_v} = \| \tilde D_u^{s-1} \phi\|_{L^2_u 
L^\infty_v}, \quad
\| \phi \|_{H^{s-1}_v L^\infty_u} = \| \tilde D_v^{s-1} \phi\|_{L^2_v L^\infty_u}.
\end{equation}
The estimate \eqref{hs-uv} implies that $\phi$ is in
$H^{s,s}$, and \eqref{bs-cons} would follow from \eqref{xd-wave}.
Conversely, when $T_0$ is small then
Theorem \ref{lwp-reg} implies that the solution $\phi$ is in
$H^{s,s}$, and the above claims follow from Sobolev embedding.

Of course, this derivation of \eqref{hs-u}- \eqref{hs-uv} is not rigorous since
the wave map equation \eqref{hardwire}, which gave \eqref{h1-u}, does not
commute with fractional integration operators as $R$ is not
constant coefficient.  However one may hope to
obtain some regularity control on $R$ and thus obtain an approximate
conservation law, using paraproduct type estimates to control the error.
It turns out that when $T_0$ is large one needs to first rescale
the solution as in \eqref{scaling-def} in order to make this approach viable.
We give rescaling precedence over differentiation, hence
$\phi^\lambda_u = (\phi^\lambda)_u$.

The rescaled versions of \eqref{hs-u}, \eqref{hs-v}, \eqref{hs-uv}
that we will rigorously prove are as follows.
\begin{theorem}\label{bs-global}  Let $3/4 < s < 1$, $C,M>0$, $T_0 \gtrsim 1$ 
be fixed, and 
suppose that the initial data to the Cauchy problem \eqref{wave-map}
satisfies \eqref{compact}, \eqref{hs-est}, \eqref{h1-finite}.  Then the 
$H^1$ solution
$\phi$ to \eqref{wave-map} satisfies the global estimates
\begin{align}
\|\phi^\lambda_u \|_{H^{s-1}_u L^\infty_v} &\leq C_1 \lambda^{\frac{1}{2}-s}
\label{bs-u}\\
\|\phi^\lambda_v \|_{H^{s-1}_v L^\infty_u} &\leq C_1 \lambda^{\frac{1}{2}-s}
\label{bs-v}\\
\|\phi^\lambda_{uv} \|_{H^{s-1}_u H^{s-1}_v} &\leq C_1 \lambda^{1-2s}
\label{bs-uv}
\end{align}
for $\lambda = C_2 T_0^{N_2}$, where $C_1, C_2, N_2 > 0$ are 
constants which do not 
depend on $T_0$ or $M$.
\end{theorem}

Note that we have the scaling relationship
\be{scaling-est}
 \lambda^{1/2} \|f\|_{H^{s-1}} \lesssim \| f^\lambda \|_{H^{s-1}} \lesssim \lambda^{\frac{3}{2}-s} \|f\|_{H^{s-1}}
\end{equation}
when $\lambda \gg 1$, and similarly for $\phi$.  Thus the estimates
\eqref{bs-u}, \eqref{bs-v}, \eqref{bs-uv} are implied by, but are weaker
than, their $\lambda=1$ counterparts \eqref{hs-u}, \eqref{hs-v}, 
\eqref{hs-uv}, especially for the
low frequency modes of $\phi$. This will be  enough to recover
polynomial growth of the $H^s$ norm, since for frequencies which are 
$\gg \lambda$ the two estimates are essentially equivalent.  

The general approach implicit in Theorem \ref{bs-global}, is motivated  
by that in \cite{bourg.2d}, where Bourgain shows
global well-posedness results (for the NLS and NLW equations) in 
spaces rougher than the energy space. Note that since $H^1$ solutions 
remain in $H^1$, the agent of blowup
in $H^s$ for $s < 1$ must be the migration of energy from high 
frequencies to low frequencies.   The bounds \eqref{bs-u}-\eqref{bs-uv}
provide control on the movement from frequencies $|\xi| >> \lambda$
to lower frequencies.  

Note however that the techniques in \cite{bourg.2d} do not
apply directly to our situation since there is no smoothing 
for the one-dimensional wave-map equation; more precisely, the
estimate
$$ \| \phi(T) - S(f,g)(T) \|_{H^1} \lesssim \|f\|_{H^s} + \|g\|_{H^{s-1}}$$
which is central to the approach in \cite{bourg.2d} does not hold for 
any $s<1$, even for short times $T$.  Our approach
relies on the very strong conservation laws in Lemma \ref{pohl_law} to
overcome this lack of smoothing.

Assuming Theorem \ref{bs-global} for the moment, let us conclude
the proof of \eqref{bs-cons}
and hence Theorem \ref{main}.  
By \eqref{xd-wave} it suffices to show that
\be{eta-claim}
\| \eta \phi \|_{H^s_u H^s_v} \lesssim T_0^N,
\end{equation}
where $\eta$ is a cutoff function adapted to the 
diamond $\{ (u,v): |u|, |v| \lesssim T_0\}$.  On the other hand, 
applying \eqref{bs-uv}
and \eqref{scaling-est}, we obtain (after expanding $\lambda$ in terms of $T_0$)
$$ \| \phi_{uv} \|_{H^{s-1}_u H^{s-1}_v} \lesssim T_0^N.$$
Since $\phi = S(f,g) + \Box^{-1} \phi_{uv}$, the claim \eqref{eta-claim} 
follows from Lemma \ref{integ} and \eqref{dx-wave}.
\section{Localized $H^s$, and one-dimensional 
paraproduct estimates}\label{wh-sec}

In this section $1 > s > 1/2$ is fixed.
In the local well-posedness theory developed in previous sections, estimates
such as $H^s H^s \subset H^s$, $H^s H^{s-1} \subset H^{s-1}$ (together
with product space analogues) were crucial.  In order to show global
well-posedness we will need to strengthen these inclusions in a number of 
ways.

Our first observation is that we may replace the space
$H^s$ by a localized variant, which we denote by $L$.
We cover the real line by finitely overlapping intervals $\{I\} = \{J\}$
of length approximately 1, and for each $I$ let $\eta_I$ be a standard 
bump function adapted to $I$ so that $\sum_{I} \eta_I \sim 1$.
\begin{definition}  If $f(u)$ is a test function, define 
the $L = L_u$ norm by
$$
\| f \|_{L_u} = \sup_I \| \tilde D_u^s (\eta_I f)\|_{L^2_u}.
$$
where $\tilde D_u$ is defined in Definition \ref{dtil}.  Similarly we define
$L = L_v$ for functions of $v$.
\end{definition}

Note that Corollary \ref{stable} implies that this definition is independent
of the exact choice of $\eta_I$.

The Sobolev spaces $H^\stil$ can be described locally as follows.

\begin{lemma}\label{local} Let $\stil$ be any real number.
If $f(u)$ is a test function, we have
\be{h-local}
 \|f\|_{H^\stil_u} \sim (\sum_I \| \eta_I f\|_{H^\stil_u}^2)^{1/2}.
\end{equation}
If $\phi(u,v)$ is a test function, 
$$ \| \phi\|_{H^\stil_u H^\stil_v} 
\sim (\sum_I \sum_J \| \eta_I(u) \eta_J(v) \phi \|_{H^\stil_{u,v}}^2)^{1/2}.$$
\end{lemma}

\begin{proof}
We prove \eqref{h-local}; the second  estimate is proven similarly.

Suppose first that $f \in H^\stil$.  By Plancherel's theorem, one can
write $f = \tilde D_u^{-\stil} F$ for some $F \in L^2$, with $\|F\|_2 \sim \|f\|_{H^\stil}$.  We may write $F = \sum_I
F_I$ where each $F_I$ is supported in $I$ and $\|F\|_2 \sim (\sum_I \|F_I\|_2^2)^{1/2}$.  Thus for any $I$
$$ \eta_I f = \sum_J \eta_I (D_u^{-\stil} F_J).$$
Since $\tilde D_u^{-\stil}$ has compactly supported kernel, the summands will
vanish unless $\dist(I,J) \lesssim 1$.  If we now invoke the triangle inequality
and discard the $\eta_I$ cutoff by Corollary \ref{stable}, we have
$$ \| \eta_I f \|_{H^\stil} \lesssim \sum_{J: \dist(I,J) \lesssim 1} 
\| \tilde D_u^{-\stil} F_J \|_{H^\stil} \sim 
\sum_{J: \dist(I,J) \lesssim 1} \|F_J\|_{2}.$$
and \eqref{h-local} follows since there are only a finite number of $J$
for each $I$.

Conversely, suppose that $(\sum_I \| \eta_I f\|_{H^\stil}^2)^{1/2}$ is finite.
We may write
$$ f = \sum_I \tilde \eta_I \eta_I f$$
for some cutoffs $\tilde \eta_I$ which are adapted to slight dilates of
$I$.  We have to estimate
$$ \|f\|_{H^\stil} \sim \| \tilde D_u^{\stil} f\|_2
= \| \sum_I \tilde D_u^{\stil} \tilde \eta_I \eta_I f \|_2.$$
Since $\tilde D_u^\stil$ has compactly supported kernel, the summands
are supported on slight dilates of $I$, and are therefore finitely overlapping.  Thus we have
$$ \|f\|_{H^\stil} \lesssim (\sum_I \| \tilde D_u^{\stil} \tilde \eta_I \eta_I f \|_2^2)^{1/2} \sim (\sum_I \| \tilde \eta_I \eta_I f\|_2^2)^{1/2},$$
and the claim then follows from Corollary \ref{stable}.
\end{proof}

There is a slight subtlety involved in product norms involving 
the space $L$.  Note that in the definition below the 
$\sup$ is inside the sum.
\begin{definition}  If $\phi(u,v)$ is a test function and $\stil 
\in \R$, we define the norms $H^\stil_u L_v$, 
$H^\stil_v L_u$, $L_u L_v = L_v L_u$ by
\bas
\| \phi \|_{H^\stil_u L_v} &=  (\sum_I \sup_J \| \eta_I(u) \eta_J(v) \phi
\|_{H^\stil_{u}H^s_v}^2)^{1/2} \\
\| \phi \|_{H^\stil_v L_u} &=   (\sum_J \sup_I \| \eta_I(u) \eta_J(v) \phi
\|_{H^s_u H^\stil_v}^2)^{1/2} \\
\| \phi \|_{L_u L_v} &=
\sup_I \sup_J \| \tilde D_u^s \tilde D_v^s (\eta_I(u) \eta_J(v) 
\phi) \|_{L^2_u L^2_v}.
\end{align*}
\end{definition}
 
We now prove some algebraic relationships between $L$ and the
Sobolev space $H^{s-1}$. The following Lemma contains  localized
variants of the embedding $H^s H^{s-1} \subset
H^{s-1}$ given in  Lemma \ref{easy-mult} above.
\begin{lemma}\label{product-lemma}  If $f(u)$ and $g(u)$ are test functions,
then
\be{h-wh}
\| fg\|_{H^{s-1}} \lesssim \|f\|_L \|g\|_{H^{s-1}}.
\end{equation}
Furthermore, if $\phi(u,v)$ and $\psi(u,v)$ are test functions, then
\ba
\| \phi \psi\|_{H^{s-1}_u H^{s-1}_v} &\lesssim \| \phi \|_{L_u L_v} \|\psi\|_{H^{s-1}_u H^{s-1}_v}
\label{hh-wwhh}\\
\| \phi \psi\|_{H^{s-1}_u H^{s-1}_v} &\lesssim \| \phi \|_{H^{s-1}_u L_v} \|\psi\|_{H^{s-1}_v L_u}
\label{hh-hwhw}
\end{align}
\end{lemma}

\begin{proof}
We prove only \eqref{h-wh}; the other two estimates follow by arguing 
similarly in both $u$ and $v$.

From \eqref{h-local} (with $\eta_I$ replaced by $\eta_I^2$) we have
$$ \|  fg\|_{H^{s-1}} \lesssim (\sum_I \| \eta_I f \eta_I g \|_{H^{s-1}}^2)^{1/2}.$$
But from Lemma \ref{easy-mult} and the definition of $L$ we have
$$  \| \eta_I f \eta_I g \|_{H^{s-1}} \lesssim
\| \eta_I f\|_{H^s} \| \eta_I g\|_{H^{s-1}} \lesssim \|f\|_L \| \eta_I g\|_{H^{s-1}}.$$
Combining this with the above estimate and using \eqref{h-local} again
one obtains \eqref{h-wh}.
\end{proof}

In the sequel we will attempt to commute integration operators such
as $\tilde D_u^{s-1}$ with identities such as \eqref{hardwire}.  In doing so
it will be natural to try to control paraproduct expressions such as
$$ \tilde D_u^{s-1}(\phi \psi) - \phi \tilde D_u^{s-1}(\psi),$$
in terms of $\phi(u,v)$, $\psi(u,v)$.  
This quantity is of comparable strength to $\tilde D_u^{s-1}(\phi \psi)$, but exhibits
cancellation when $\phi$ is constant or slowly varying.

In the next section we will need to
 estimate the above quantity in $L^2_u H^{s-1}_v$.
If one ignored the cancellation and used the triangle inequality, one
obtains
$$ \| \tilde D_u^{s-1}(\phi \psi) - \phi \tilde D_u^{s-1} \psi\|_{L^2_u H^{s-1}_v} \lesssim
\| \phi \psi\|_{H^{s-1}_u H^{s-1}_v} + \| \phi \tilde D_u^{s-1} \psi\|_{L^2_u H^{s-1}_v}.$$
Combining this with Lemma \ref{product-h1} one obtains the bound
\begin{equation}
\label{see.me}
\| \tilde D_u^{s-1}(\phi \psi) - \phi \tilde D_u^{s-1} \psi\|_{L^2_u H^{s-1}_
v} \lesssim
 \| \phi \|_{H^s_u H^{s-1}_v} \| \psi\|_{H^{s-1}_u H^s_v}.
\end{equation}
By the previous discussion, we may improve this estimate by
localizing the $H^s_v$ norm to an $L_v$ norm.  We could also improve
the $H^s_u$ norm in this manner, but we will instead pursue a different 
improvement
which tries to take advantage of the cancellation if $\phi$
has low frequency.  In fact, we have

\begin{lemma}\label{exotic} 
If $\phi(u,v)$ and $\psi(u,v)$ are test functions, then
\be{exotic-est}
\| \tilde D_u^{s-1}(\phi \psi) - \phi \tilde D_u^{s-1}(\psi) \|_{L^2_u H^{s-1}_v} \lesssim 
\| \phi_u \|_{H^{s-1}_u H^{s-1}_v} \| \psi \|_{H^{s-1}_u L_v}.
\end{equation}
\end{lemma}

Note that $\|f_u \|_{H^{s-1}}$ is essentially the same as $\|f\|_{H^s}$
when $f$ consists of high frequencies, but is somewhat smaller
for low frequencies, in accordance with the previous heuristics
concerning the cancellation.

\begin{proof}  The first step is to replace the $L_v$ norm with the
stronger $H^s_v$ norm.
Let $\tilde{\eta}_J$ be a cut-off function which is one on the support
of $\eta_J$, as in the proof of Lemma \ref{local}.   
It suffices to prove the estimate
\begin{align}
\|\eta_J(v) [\tilde D_u^{s-1}(\phi \psi) - 
\phi \tilde D_u^{s-1}(\psi)] \|_{L^2_u H^{s-1}_v} & \lesssim 
\| \tilde{\eta}_J(v) \phi_u \|_{H^{s-1}_u H^{s-1}_v} \| 
\eta_J(v) \psi \|_{H^{s-1}_u H^s_v}. \label{middy}
\end{align}
uniformly in $J$, since \eqref{exotic-est} can be recovered by square-summing
\eqref{middy} in $J$, using the compact support of the kernel of 
$\tilde D_u^{s-1}$, and applying Lemma \ref{local}.  
By the support properties of $\eta, \tilde{\eta}$, we may rewrite
\eqref{middy} as
$$
\|[\tilde D_u^{s-1}(\tilde\eta_J(v) \phi \eta_J(v) \psi) - 
\tilde \eta_J(v) \phi \tilde D_u^{s-1}(\eta_J(v) \psi)] 
\|_{L^2_u H^{s-1}_v} \lesssim 
\| \tilde{\eta}_J(v) \phi_u \|_{H^{s-1}_u H^{s-1}_v} \| 
\eta_J(v) \psi \|_{H^{s-1}_u H^s_v}. 
$$
Replacing $\tilde{\eta}_J \phi$ with 
$\phi$ and $\eta_J \psi$ with $\psi$, it
suffices to show 
\be{exotic-local}
\| \tilde D_u^{s-1}(\phi \psi) - \phi \tilde D_u^{s-1}(\psi) \|_{L^2_u H^{s-1}_v} \lesssim 
\| \phi_u \|_{H^{s-1}_u H^{s-1}_v} \| \psi \|_{H^{s-1}_u H^s_v}.
\end{equation}
for arbitrary test functions $\phi$, $\psi$.  (This estimate
should be compared with \eqref{exotic-est}).

By Plancherel's theorem, the left-hand side is equal to the $L^2_\mu L^2_\nu$
norm of
\be{quant}
 C \int_{\mu_1 + \mu_2 = \mu} \int_{\nu_1 + \nu_2 = \nu}
[m_{s-1}(\mu_1 + \mu_2) - m_{s-1}(\mu_2)] m_{s-1}(\nu_1+\nu_2) 
\hat\phi(\mu_1,\nu_1) \hat\psi(\mu_2,\nu_2)\ d\mu_1 d\nu_1.
\end{equation}
Define $\Phi$, $\Psi$ by $\hat \Phi(\mu_1,\nu_1) = |\mu_1| \langle \mu_1\rangle^{-1} |\hat \phi(\mu_1,\nu_1)|$, $\hat \Psi(\mu_2,\nu_2) = |\hat \psi(\mu_2,\nu_2)|$; note that
\be{Phi}
 \| \Phi \|_{H^s_u H^{s-1}_v} \sim \| \phi_u \|_{H^{s-1}_u H^{s-1}_v},
\quad \| \Psi \|_{H^{s-1}_u H^s_v} \sim \| \psi \|_{H^{s-1}_u H^s_v}.
\end{equation}
The quantity \eqref{quant} is majorized by
\be{quant2}
 \int_{\mu_1 + \mu_2 = \mu} \int_{\nu_1 + \nu_2 = \nu}
\frac{|m_{s-1}(\mu_1 + \mu_2) - m_{s-1}(\mu_2)|}{|\mu_1|} \langle \mu_1 \rangle \langle \nu_1 + \nu_2\rangle^{s-1}
\hat \Phi(\mu_1,\nu_1) \hat \Psi(\mu_2,\nu_2)\ d\mu_1 d\nu_1.
\end{equation}
When $|\mu_1| \gtrsim 1$ we have
$$ \frac{|m_{s-1}(\mu_1 + \mu_2) - m_{s-1}(\mu_2)|}{|\mu_1|} \langle \mu_1\rangle \lesssim \langle \mu_2 \rangle^{s-1} + \langle \mu_1 + \mu_2 \rangle^{s-1},
$$
while when $|\mu_1| \lesssim 1$ the mean-value theorem gives
$$ \frac{|m_{s-1}(\mu_1 + \mu_2) - m_{s-1}(\mu_2)|}{|\mu_1|} \langle \mu_1\rangle \lesssim \langle \mu_2 \rangle^{s-2}.
$$
Thus in either case we have
$$ \frac{|m_{s-1}(\mu_1 + \mu_2) - m_{s-1}(\mu_2)|}{|\mu_1|} \langle \mu_1
\rangle \lesssim \langle \mu_2 \rangle^{s-1} + \langle \mu_1 + \mu_2 \rangle^{s-1}.
$$
Inserting this into \eqref{quant2} and using Plancherel's theorem, we see
that the $L^2_\mu L^2_\nu$ norm of \eqref{quant2} is majorized by
$$ \| D_v^{s-1} (\Phi D_u^{s-1} \Psi) \|_{L^2_u L^2_v}
+ \| D_u^{s-1} D_v^{s-1} (\Phi \Psi) \|_{L^2_u L^2_v} 
= \| \Phi D_u^{s-1} \Psi \|_{H^{s-1}_u L^2_v} + \| \Phi \Psi \|_{H^{s-1}_u
H^{s-1}_v}.$$
By Lemma \ref{product-h1} this is majorized by
$$ \| \Phi\|_{H^s_u H^{s-1}_v} \| \Psi \|_{H^{s-1}_u H^s_v},$$
and the claim now follows from \eqref{Phi}.
\end{proof}

To close this section we give some elementary estimates which connect
the $L$ space to $H^{s-1}$ and the $L^\infty$ norm; this
will allow us to translate the estimates in Theorem \ref{bs-global}
to ones involving $L$.

\begin{lemma}\label{l-conv}  If $f(u)$ is a test function, then
\be{l-ih}
\| f\|_{L_u} \lesssim \|f\|_\infty + \|f_u\|_{H^{s-1}_u}.
\end{equation}
If $\phi(u,v)$ is a test function, then
\bas
 \| \phi\|_{H^{s-1}_v L_u} &\lesssim \|\phi\|_{H^{s-1}_v L^\infty_u} + 
\|\phi_u\|_{H^{s-1}_v H^{s-1}_u}\\
 \| \phi\|_{H^{s-1}_u L_v} &\lesssim \|\phi\|_{H^{s-1}_u L^\infty_v} + 
\|\phi_v\|_{H^{s-1}_u H^{s-1}_v}\\
 \| \phi\|_{L_u L_v} &\lesssim \|\phi\|_{L^\infty_u L^\infty_v}
+ \| \phi_u \|_{H^{s-1}_u L^\infty_v} + \| \phi_v \|_{H^{s-1}_v L^\infty_u}
+ \| \phi_u \phi_v \|_{H^{s-1}_u H^{s-1}_v}
\end{align*}
where the norms $H^{s-1}_u L^\infty_v$, $H^{s-1}_v L^\infty_u$ were
defined in \eqref{hl-def}.
\end{lemma}

\begin{proof}
We prove only \eqref{l-ih}; the other estimates follow by applying
a similar argument applied to both variables at once. 

It suffices to show that
$$ \| \eta_I f\|_{H^s} \lesssim \|f\|_\infty + \|f_u\|_{H^{s-1}_u}.$$
We may partition frequency space and divide $f$ into a piece with
frequency support on $|\mu| \lesssim 1$, and a piece with frequency
support on $|\mu| \gtrsim 1$.  To handle the first piece we use the
estimate
$$  \| \eta_I f\|_{H^s} \lesssim \| \eta_I f\|_{H^N} \lesssim \| f\|_{C^N}
\lesssim \|f\|_\infty$$
for some large integer $N$, where the last inequality follows from the
frequency support hypothesis.

To handle the second piece we use Lemma \ref{easy-mult} to obtain
$$  \| \eta_I f\|_{H^s} \lesssim \| f\|_{H^s} \sim \| f_u \|_{H^{s-1}}$$
where the last inequality follows from the frequency support hypothesis.
\end{proof}

\begin{lemma}\label{compact-est} Suppose that $f(u)$ is a test function supported
on an interval $I$ of length $\gtrsim 1$.  Then
$$ \| f\|_{L^\infty} \lesssim \|f\|_{L^\infty(I^\prime)} + 
|I|^{1/2} \|f_u\|_{H^{s-1}_u},$$
where $I^\prime$ is any nonempty subinterval of $I$.
\end{lemma}

\begin{proof}
It suffices to show that
$$ |f(u) - f(u_0)| \lesssim |I|^{1/2} \|f_u\|_{H^{s-1}_u}$$
whenever $u \in I$, $u_0 \in I^\prime$.  But by the fundamental theorem
of calculus the left hand side is majorized by
$$ |\langle f_u, \chi_{[u_0,u]} \rangle| \lesssim \| f_u \|_{H^{s-1}}
\| \chi_{[u_0,u]} \|_{H^{1-s}},$$
and the result follows from the hypothesis $s > 1/2$ and the
inequality $\| \chi_{[u_0,u]} \|_{H^{1-s}} \lesssim \langle u_0 - 
u \rangle^{1/2} \lesssim |I|^{1/2}$.
When $u_0 = -1$, $u=1$ this inequality follows from direct computation,
and the general case follows by rescaling and translation invariance.
\end{proof}

\section{Proof of Theorem \ref{bs-global}}\label{proof-sec}

Fix $C_0$, $T_0$, $M$, $3/4 < s < 1$.  We will let $C_1$ be a large
constant to be chosen later, and $C_2$, $N_2$ to be large constants
depending on $C_1$, also to be chosen later.  In particular, $\lambda$
is also fixed.  The quantities $N$, $C$, $\eps$ and the implicit constants
in the estimates will vary from line to line, but will not depend on $C_2$.

We shall use the continuity method.  Let $B$ denote the set
$$ B = \{ (f,g): \hbox{\eqref{compact}, \eqref{hs-est}, 
\eqref{h1-finite} hold.} \}.$$
We give $B$ the induced topology from $H^1 \times L^2$.
Consider the subset of $B$
$$ E = \{ (f,g) \in B: \hbox{\eqref{bs-u}, \eqref{bs-v}, \eqref{bs-uv} hold.}
\}.$$
We wish to show that $E=B$.  To this end we introduce the
weaker versions of \eqref{bs-u}, \eqref{bs-v}, \eqref{bs-uv}
\begin{align}
\|\phi^\lambda_u \|_{H^{s-1}_u L^\infty_v} &\lesssim \lambda^{\frac{1}{2}-s+\eps}
\label{bs-u2}\\
\|\phi^\lambda_v \|_{H^{s-1}_v L^\infty_u} &\lesssim \lambda^{\frac{1}{2}-s+\eps}
\label{bs-v2}\\
\|\phi^\lambda_{uv} \|_{H^{s-1}_u H^{s-1}_v} &\lesssim \lambda^{1-2s+2\eps}.
\label{bs-uv2}
\end{align}
and define the subset of $B$
$$ \tilde E = \{ (f,g) \in B: \hbox{\eqref{bs-u2}, \eqref{bs-v2}, \eqref{bs-uv2} hold}\}.$$
Clearly $E \subset \tilde E$ if $C_2$, $N_2$ are sufficiently large.  
Furthermore, we claim the following:
\begin{itemize}
\item  If $C_2$, $N_2$ are sufficiently large, then there exists
an $\epsilon_M > 0$ which can depend on $T_0$, $M$, $\lambda$ 
such that the following
holds: If $(f,g) \in E$ and $(\tilde f, \tilde g)$ is within $\epsilon_M$
of $(f,g)$ in $H^1 \times L^2$ norm, then $(\tilde f, \tilde g)$ is
in $\tilde E$.
\item  If $(f,g)$ is in $\tilde E$, then $(f,g)$ is in $E$.
\end{itemize}
Combining these two statements we see that $E$ is both open and closed
in $H^1 \times L^2$.  Since $E$ contains the origin and
$B$ is connected, we will be done.

To prove the first claim, we first observe that 
\eqref{bs-u}, \eqref{bs-v}, \eqref{bs-uv}
are trivial to verify outside of the diamond $\{|u|, |v| \lesssim T_0\}$,
by Corollary \ref{scatter}.  Thus we may restrict our attention to the diamond,
which is a compact set.  

From Theorem \ref{lwp-reg} we see that the $H^{1,1}$ norm of $\phi^\lambda$
on the diamond depends in a Lipschitz manner on the $H^1 \times L^2$ norm of 
the data (with a large Lipschitz constant depending on $M$, $T_0$, $\lambda$).
Since the $H^{1,1}$ norm controls the norms present in the definition
of $E$, $\tilde E$ by Sobolev embedding, the claim follows by elementary
topology.

The remainder of this section is devoted to proving the second claim.
Accordingly, we fix $(f,g) \in B$, assume that \eqref{bs-u2}, \eqref{bs-v2}, 
\eqref{bs-uv2} hold,
and try to prove \eqref{bs-u}, \eqref{bs-v}, and \eqref{bs-uv}.

Since $\phi^\lambda$ stays on the sphere, we have
\be{bs-infty}
\|\phi^\lambda \|_{L^\infty_u L^\infty_v} \lesssim 1.
\end{equation}
Since $\lambda$ is large and $\frac{1}{2} - s + \eps < 0$ for $\eps$ 
sufficiently small, we can use Lemma \ref{l-conv} and 
\eqref{bs-infty}, \eqref{bs-u2}, \eqref{bs-v2}, \eqref{bs-uv2} to obtain
estimates involving the space $L$.  More precisely, we have
\ba
\| \phi^\lambda \|_{L_u L_v} &\lesssim 1 \label{ww-est}\\
\| \phi^\lambda_u \|_{H^{s-1}_u L_v} &\lesssim \lambda^{\frac{1}{2}-s+\eps} \label{hw-est}\\
\| \phi^\lambda_v \|_{H^{s-1}_v L_u} &\lesssim \lambda^{\frac{1}{2}-s+\eps} \label{wh-est}\\
\| \phi^\lambda_{uv} \|_{H^{s-1}_u H^{s-1}_v} &\lesssim \lambda^{1-2s+2\eps} \label{hh-est}.
\end{align}
We now show \eqref{bs-u}.  
This is the same (if $C_1$ is chosen sufficiently large) as
\be{start}
 \| \tilde D_u^{s-1} \phi^\lambda_u \|_{L^2_u L^\infty_v} 
\lesssim \lambda^{\frac{1}{2}-s}.
\end{equation}

When $s=1$ this was proven in Section \ref{global-h1-sec} by the computation
\eqref{asym}.  The argument here will be an adaptation of this computation.

We first prove \eqref{start} for short times $|t| \lesssim \eps \lambda$,
i.e. we show
$$ \| \chi^\lambda \tilde D_u^{s-1} \phi^\lambda_u \|_{L^2_u L^\infty_v}
\lesssim \lambda^{\frac{1}{2}-s},$$
where $\chi$ is a cutoff which equals one the slab $|t| \lesssim \eps$,
and vanishes on a dilate of this slab.
Since $\tilde D_u^{s-1}$ has compactly supported kernel,
we may write this as
\be{first-part}
 \| \chi^\lambda \tilde D_u^{s-1} (\tilde \chi \phi)^\lambda_u \|_{L^2_u L^\infty_v}
\end{equation}
where $\tilde \chi$ equals 1 on a dilate on the support of $\chi$, 
and vanishes outside of an even larger dilate.  
Discarding
the $\chi^\lambda$ term and using \eqref{scaling-est}, we reduce
ourselves to showing that
$$ \| (\tilde \chi \phi)_u \|_{H^{s-1}_u L^\infty_v} \lesssim 1.$$
However, for short times $|t| \lesssim \eps$ Theorem \ref{lwp-reg} 
applies, and we have
$$ \| \tilde \chi \phi\|_{H^s_u H^s_v} \lesssim 1.$$
The claim then follows by Sobolev embedding.

We now prove the full estimate \eqref{start}.  By squaring, we obtain
\be{start2}
\| |\tilde D_u^{s-1} \phi^\lambda_u|^2 \|_{L^1_u
L^\infty_v} \lesssim \lambda^{1-2s}.
\end{equation}
Since this estimate was just proven for short times, we may invoke
Lemma \ref{compact-est}, and reduce ourselves to showing that
$$ \| \partial_v |\tilde D_u^{s-1} \phi^\lambda_u|^2
\|_{L^1_u H^{s-1}_v} \lesssim (\lambda T_0)^{-1/2} \lambda^{1-2s}.$$
By evaluating the $v$ derivative and using \eqref{hardwire}, 
it thus suffices to show that
\be{easy-est}
 \| (\tilde D_u^{s-1} \phi^\lambda_u)^t (\tilde D_u^{s-1} R \phi^\lambda_u) 
\|_{L^1_u H^{s-1}_v} \lesssim (\lambda T_0)^{-1/2} \lambda^{1-2s}
\end{equation}
Since $R$ is anti-symmetric,
$$  (\tilde D_u^{s-1} \phi^\lambda_u)^t R (\tilde D_u^{s-1} \phi^\lambda_u) = 0.$$
Thus it suffices to show that
$$ \| (\tilde D_u^{s-1} \phi^\lambda_u)^t [(\tilde D_u^{s-1} R \phi^\lambda_u)
- R (\tilde D_u^{s-1} \phi^\lambda_u)] \|_{L^1_u H^{s-1}_v} \lesssim (\lambda T_0)^{-1/2} 
\lambda^{1-2s}.$$
However, from \eqref{bs-u2} we have
$$ \| \tilde D_u^{s-1} \phi^{\lambda}_u \|_{L^2_u L_v} \lesssim
\lambda^{\frac{1}{2} -s + \eps}.$$
Thus by \eqref{h-wh} it suffices to show that
$$ \| \tilde D_u^{s-1} (R \phi^\lambda_{u}) - R \tilde D_u^{s-1} 
\phi^\lambda_{u} 
\|_{L^2_u H^{s-1}_v} 
\lesssim (\lambda T_0)^{-1/2} \lambda^{\frac{1}{2} -s -\eps}.$$
By Lemma \ref{exotic} this reduces to
$$ \| R_u \|_{H^{s-1}_v H^{s-1}_u} \| \phi^\lambda_u \|_{H^{s-1}_v L_u} \lesssim (\lambda T_0)^{-1/2} \lambda^{\frac{1}{2} -s -\eps}.$$
By \eqref{bs-u2} again, it thus suffices to show that
$$ \| R_u \|_{H^{s-1}_v H^{s-1}_u} \lesssim (\lambda T_0)^{-1/2} \lambda^{-2\eps}.$$
By expanding out $R_u$, we need only show
$$ \| (\phi^\lambda)^t \phi^\lambda_{uv} \|_{H^{s-1}_v H^{s-1}_u} + \| (\phi^\lambda_u)^t \phi^\lambda_v \|_{H^{s-1}_v H^{s-1}_u}
 \lesssim (\lambda T_0)^{-1/2} \lambda^{-2\eps}.$$
But by \eqref{hh-wwhh}, \eqref{hh-hwhw}, the left-hand side of this
is majorized by
$$ \| \phi^\lambda \|_{L_v L_u} \| \phi^\lambda_{uv} \|_{H^{s-1}_v H^{s-1}_u} 
+ \|  \phi^\lambda_u \|_{H^{s-1}_u L_v} \|  \phi^\lambda_v \|_{H^{s-1}_v L_u},$$
which is bounded by $\lambda^{1-2s + 2\eps}$ by \eqref{ww-est}, \eqref{wh-est},
\eqref{hw-est}, \eqref{hh-est}.  Since $s > 3/4$, the claim is thus proven if
$\eps$ is sufficiently small and $C_2$, $N_2$ are sufficiently large.
This concludes the proof of \eqref{bs-u}. Note that if one used \eqref{see.me}
instead of Lemma \ref{exotic} then we'd 
need $s>1$ instead of $s > \frac{3}{4}$. 

The proof of \eqref{bs-v} is similar, so we turn to \eqref{bs-uv}.  It suffices
to show that (if $C_1$ is sufficiently large)
$$
\| \phi^\lambda_{uv} \|_{H^{s-1}_u H^{s-1}_v} \lesssim \lambda^{1-2s}.
$$
From \eqref{wave-map} the left-hand side is
$$
\| \phi^\lambda (\phi^\lambda_u)^t \phi^\lambda_v \|_{H^{s-1}_u H^{s-1}_v}.$$
But by \eqref{hh-wwhh}, \eqref{hh-hwhw}, this is majorized by
$$ \| \phi^\lambda \|_{L_u L_v} \|\phi^\lambda_u\|_{H^{s-1}_u L_v} \| \phi^\lambda_v\|_{H^{s-1}_v L_u}.$$
Thus by \eqref{ww-est} it suffices to show that
$$ \|\phi^\lambda_u\|_{H^{s-1}_u L_v} \lesssim \lambda^{\frac{1}{2}-s}$$
and similarly for $\phi^\lambda_v$.  But this follows from the definition
of $H^{s-1}_u L_v$, the estimate \eqref{bs-u} just proved, and
\eqref{bs-uv2} (if $\eps$ is sufficiently small).
\endprf

\section{Negative results}\label{negative}

In this section we give some rather simple negative results regarding
ill-posedness of the wave map equation and similar equations.  

The nonlinearity in the wave map equation \eqref{wavemap.null} 
contains the null form
$Q_0^{\alpha\beta}(\phi,\phi) \equiv \phi^\alpha_u \phi^\beta_v$.
That the quadratic form $Q_0$ has this null structure is
important for low regularity well-posedness, as the following 
simple example shows.  (See \cite{lindblad}, \cite{kman.mach.smoothing}
for a similar situation in dimension $n=3$.) 

\begin{proposition}  The scalar equation
$$ \Box \phi = |\phi_v|^2$$
is locally well-posed in $H^s$ if and only if $s > 3/2$.
\end{proposition}

\begin{proof}  By making the substitution $\psi = \phi_v$,
it suffices to show that the equation
$$ \psi_u = |\psi|^2$$
is locally well posed in $H^s$ if and only if $s > 1/2$.  
But from the explicit solution
$$ \psi(x,t) = \frac{\psi(x-t,0)}{1 - \psi(x-t,0) t},$$
we see that $\psi$ only stays regular for a non-zero time when the 
initial date $\psi(\cdot, 0)$ is bounded.  This is only guaranteed 
when $s > 1/2$,
hence the result.
\end{proof}

With the null form structure, one can do much better, as the following
example of Nirenberg shows.  
\begin{proposition} \cite{kman.mach.nullform}  The scalar equation
\begin{align*}
\Box \phi & = \phi_u \phi_v \\
\phi(x,0) \; = \; f\;  &  \; \phi_t(x,0) \;  = \;  g
\end{align*}
is locally well-posed in $H^s$ if and only if $s > 1/2$.
\end{proposition}

\begin{proof} For completeness, we sketch the argument given 
in \cite{kman.mach.nullform} here.
Take data $f = 0$, $g \in H^{s -1}$. By making the 
substitution $\psi = 1 - e^\phi$, 
it suffices to show that the solution to
$$ \Box \psi = 0$$
remains in $H^s$ and satisfies $\| \psi(t,x) \|_{L^{\infty}(\rr)} < 1$
for a non-zero amount of time.  This
is true for $s>1/2$ since $\psi$ is then uniformly continuous
by Sobolev embedding.  For $s \leq 1/2$ it is easy to construct
discontinuous $\psi$ which becomes large instantaneously.
\end{proof}

By Theorem \ref{lwp-reg} and previously mentioned  work in higher dimensions,
one has local well-posedness in $H^s(\rr^n)$ for wave maps from 
$\rr^{n+1}$ when $s > n/2$.  For
$s < n/2$ the problem is supercritical and well-posedness seems very unlikely. 
(For $n \geq 3$, the supercritical wave map problem is ill-posed for
certain manifolds: see Shatah, Shatah-Zadeh 
\cite{shatah, shatah.zadeh.blow}.)  The critical case $s=n/2$ 
seems very subtle, as the following example demonstrates.
\begin{proposition}  \label{nosmooth} If the target manifold 
of the wave map \eqref{mapping}
is the sphere $S^{m-1}$,  $m \geq 2$, then there exist coordinates
for which the solution operator to the wave map equation \eqref{wavemap} is not twice differentiable on 
the data space $H^{\frac{n}{2}}(\rr^n)$.   
\end{proposition}

\noindent{\bf Remark}:  In particular, Proposition  \ref{nosmooth} shows
 that the solution operator in $\rr^{2+1}$ does not depend smoothly on the 
data in the energy norm, and so there exists a coordinate system on the target
manifold in which  
one cannot prove a critical $H^1$ result by the usual Picard 
iteration argument.  The proposition also holds if 
the inhomogeneous norm $H^{n/2}$ is replaced by the homogeneous 
version $\dot H^{n/2}$.

\begin{proof} It is well known (e.g. \cite{sideris.harmonic}) 
that composition of a solution of the free wave equation with 
a geodesic yields a wave map.  Let  $\psi: \rr^n \rightarrow \rr$ 
satisfy the free wave
equation with data 
\begin{align}
(\psi, \psi_t) & = (0,\epsilon g) \label{first}.
\end{align}
where $g \in H^{\frac{n}{2} - 1}$ and $\epsilon \in \rr$. 
The mapping  $\rr \rightarrow S^1 \subset S^{m-1}$ given
by $x \rightarrow e^{ix}$ is a geodesic, hence  the function
\begin{align}
\phi_{\epsilon}(x,t) & = e^{i \psi(x,t)} \; = \;  e^{i \epsilon 
\frac{\sin( \sqrt{- \Delta} t )}{\sqrt{- \Delta}} g} \label{phi}
\end{align} 
is a solution of the wave-map system \eqref{wavemapsphere} with
initial velocity $i \epsilon g$.
 
Assume for contradiction that the mapping taking initial velocity to
the solution  \eqref{phi} at time $1$
\begin{align*}
S: i\epsilon g & \rightarrow \phi_{\epsilon}(x,1) 
\end{align*}
is  twice differentiable as a mapping $S:  H^{\frac{n}{2} - 1}(\rr^n)
\rightarrow    H^{\frac{n}{2}}(\rr^n)$.  This implies
\begin{align}
\frac{d^2}{d \epsilon^2}_{\big|_{\epsilon = 0}} 
\phi_{\epsilon}(x,1) & \in H^{\frac{n}{2}}(\rr^n).
\label{double-d}
\end{align}
Hence
\begin{align}
\left(\frac{\sin(\sqrt{- \Delta})}{\sqrt{- \Delta}} g \right)^2 & 
\in H^{\frac{n}{2}}(\rr^n) 
\label{contradict.me}
\end{align}
for all $g \in H^{\frac{n}{2} - 1}(\rr^n)$. 

Consider the preliminary function  $G(x)$ defined by
\begin{align*}
\hat{G}(\xi) & = \frac{1}{ \log^{\frac{1}{2} +}(\langle \xi \rangle) 
\cdot \langle \xi \rangle ^{n-1} }.
\end{align*}
One easily verifies $G \in H^{\frac{n}{2} -1}(\rr^n)$. If we set 
$g(x) \, = \, \sin(\sqrt{- \Delta}) G(x)$ then clearly we also
have  $g \in  H^{\frac{n}{2} -1}(\rr^n)$.  A straightforward computation
gives
\begin{align*}
\| \left(\frac{\sin(\sqrt{- \Delta})}{\sqrt{- \Delta}}  
g \right)^2   
\|^2_{H^{\frac{n}{2}}} & = 
\int
\langle \xi \rangle^n \big( 
\int\frac{\sin^2(|\xi - \eta|)}{|\xi - \eta|} G(\xi - \eta) 
\frac{\sin^2(|\eta|)}{|\eta|}G(\eta) d \eta \big)^2 d \xi
\\
& \gtrsim \int \langle \xi \rangle^n
 \left( \int 
\frac{1}{|\xi - \eta|}
G(\xi - \eta) \frac{1}{|\eta|}G(\eta) d \eta \right)^2 d \xi \\
& \gtrsim \int_{|\xi| \gtrsim 1} \langle \xi \rangle^n
 \big( \int_{1 \lesssim |\eta| \lesssim |\xi|}
\frac{1}{|\xi - \eta|}
G(\xi - \eta) \frac{1}{|\eta|}G(\eta) d \eta \big)^2 d \xi \\
& \gtrsim \int_{|\xi| \gtrsim 1} \langle \xi \rangle^n \big(
\frac{1}{|\xi|} G(\xi) \int_{1}^{|\xi|}
\frac{1}{\rho}G(\rho) \rho^{n-1} d\rho \big)^2\ d\xi\\
& \gtrsim \int_{|\xi| \gtrsim 1}  
 |\xi |^n
\left( \frac{1}{|\xi|} \frac{1}{\log^{\frac{1}{2}+}(|\xi|) \cdot |\xi|^{n-1}} 
\log^{\frac{1}{2} -}(|\xi|) \right)^2 d \xi \\
& = \int_1^{\infty} \frac{1}{\rho^n \log^{0+}(\rho)} \rho^{n-1} d \rho \\
& = \infty
\end{align*}
which contradicts \eqref{contradict.me} as desired. 
\end{proof}

One can be much more precise on the nature of the solution operator.
For instance, when $n =1$, the operator is continuous but not 
uniformly continuous
on $\dot H^{1/2}$, and is neither Lipschitz nor everywhere differentiable.
Further ill-posedness results in this direction are in \cite{tao:ill}.

This example points out that the choice of coordinates or
frames on the target manifold is important.  For instance, if 
one uses intrinsic arclength coordinates on $S^1$ rather than 
extrinsic complex coordinates,
then the wave map equation becomes the free wave equation, which is of course
analytically well-posed in virtually any data space.
These conclusions are consistent with \cite{helein,grillakis.equivariant, 
christ.spherical.wave, shatah.zadeh.blow}, which work in special 
coordinate systems.  

We conclude this section with a negative scattering result.  Scattering
appears unlikely for the one-dimensional wave map, since there is no obvious
decay in the equation, and furthermore the solution stays on a manifold $\m$
while free solutions almost never do.
The following result reinforces these heuristics, at least for data 
which does not have a
conditionally integrable velocity.  In the converse direction, if
the data is compactly supported, scattering was shown in  
Corollary \ref{scatter}, and we'll show in Section \ref{l11} that one
also has scattering 
when the velocity and the derivative of the position are 
absolutely integrable.
\begin{proposition}  If the target manifold is the unit circle $S^1$
in the complex plane, the initial position is $f \equiv 1$, and
the initial velocity $ig$ is smooth, then the solution asymptotically
approaches a free solution in $\dot{H}^1(\rr)$ if and only if 
the limits $G_\pm = \lim_{x \to \pm \infty} G(x)$ exist, where $G$ 
is a primitive of $g$.  In particular, there is data in $\dot H^1
\times L^2$
which does not scatter.
\end{proposition}

\begin{proof} 
We have the explicit solution
\begin{align}
 \phi(x,t) & = e^{iG(x+t)/2} e^{-iG(x-t)/2}. \label{exacto}
\end{align}
where $G(x) \, = \, \int_0^{x} g(\lambda) d \lambda.$
If the limits $G_\pm$ exist, then it's easy to see $\phi$ approaches
the free solution
$$ \psi(x,t) = e^{iG(x+t)/2} e^{-iG_-/2} + e^{iG_+/2} e^{-iG(x-t)/2}
- e^{iG_+/2} e^{-iG_-/2}$$
in $\dot{H}^1(\rr)$ as $t \to \infty$.  For example,
\begin{align*}
\| \partial_u \left(\phi(u,v) - \psi(u,v) \right) \|_{L^2} & = 
\|g(u) e^{i G(u) /2 } \left( e^{- iG(v)/2} - e^{- iG_{-}/2} \right) \|_{L^2} \\
& \rightarrow 0 
\end{align*}
by dominated convergence. 

Suppose conversely that the solution \eqref{exacto}
approaches a free solution in $\dot{H}^1 \times L^2$:
\begin{align*} \| \nabla_{x,t} ( \phi(x,t) - f_+(x+t) - f_-(x-t) )
\|_{L^2} &  \rightarrow
 0 \quad  \text{ as } \;  t \to \infty. 
\end{align*}
The convergence of the $u$ derivative in $L^2$ gives,
\begin{align}
\| \frac{i g(u)}{2} e^{iG(u)/2} e^{-iG(u-2t)/2} - f_+^{\prime}(u) 
\|_{L^2} \to 0 \hbox{ as } t & \to \infty  \label{there.boy}
\end{align}
where we've changed variables $x = u - t$.
Suppose we restrict the $u$ integration in \eqref{there.boy} to a compact 
set on which $g(u)$ is non-zero.
Then the above convergence is only possible if $G(u-2t)$
converges as $t \to \infty$, which means that $G_-$ must exist.  A
similar argument shows that $G_+$ must also exist.
\end{proof}

\section{Global existence and scattering for large data in the critical space $L^{1,1}$}\label{l11}

We give here an elementary proof of Theorem \ref{l11-main},
which gives global existence and scattering for arbitrary target
manifolds in the critical (i.e. scale-invariant) data space $\dd = L^{1,1}$ 
defined by \eqref{l.11.def}.

As with other critical global existence results (e.g. 
\cite{grillakis.semilinear}), the proof follows a familiar
pattern:
\begin{itemize}
\item Prove global well-posedness for small $L^{1,1}$ data.
\item Extend to global well-posedness for large $L^{1,1}$ data
by a nonconcentration argument.
\end{itemize}
Scattering will be obtained by a  conformal compactification argument
in Section \ref{l11.scatter} below.
\subsection{Global existence for small data}\label{small-data}
Suppose the initial data $f,g$ is small in $L^{1,1}$ when 
measured in a single coordinate chart.

We apply Lemma \ref{black-box}
with $T_1=T_2=\infty$ and the space $X$ given by
\begin{align}
\Norm{\phi}_{X} & \equiv 
\Norm{\partial_u \partial_v \phi}_{L^1_{u,v}} +  \Norm{\phi(0)}_\dd.
\label{Xdef.L11}
\end{align}
We need to check that \eqref{data-X}, \eqref{X-data}, \eqref{algebra-property},
and \eqref{box-property} hold.

The preliminary bounds
\begin{equation}
\label{mixed.est}
\begin{split}
\Norm{\partial_u \phi}_{L^1_{u}L^{\infty}_v}
& \lesssim \Norm{\phi}_{X}\\
 \Norm{\partial_v \phi}_{L^1_{v}L^{\infty}_u}  &\lesssim \Norm{\phi}_{X}
\end{split}
\end{equation}
follow from the fundamental theorem of calculus.  For instance,
we have
\begin{align*}
\Norm{\partial_u \phi}_{L^1_{u}L^{\infty}_v} & =
\int_{-\infty}^{\infty}  \sup_v |\partial_u \phi(u,v)| du \\
& \leq \int_{-\infty}^{\infty} \left(
 \int |\partial_u \partial_v \phi(u, v')| dv' + |\partial_u\phi(u,u)| \right) du \\ 
& \leq \Norm{\phi}_X.
\end{align*}

The property \eqref{data-X} is trivial from the definition of $S(f)$, so we turn
to \eqref{X-data}.  
Using \eqref{mixed.est},
\begin{align*}
\Norm{\partial_x \phi}_{L^1_x}(T) + \Norm{\partial_t \phi}_{L^1_x}(T)
& \lesssim \Norm{\partial_u \phi (u, u - 2T)}_{L^1_u} 
+ \Norm{\partial_v \phi (v + 2T,v)}_{L^1_v} \\
& \leq \Norm{\partial_u \phi}_{L^1_u L^{\infty}_v} + \Norm{\partial_v \phi}_{L^1_v L^{\infty}_u}
\\
& \lesssim \Norm{\phi}_{X}. 
\end{align*}
To finish the proof of \eqref{X-data} we need to bound the term $\Norm{\phi(\cdot, T)}_{L^{\infty}(\rr)}$
in the norm $\Norm{\phi(\cdot, T)}_{\dd}$.  
We have
\begin{align}
\int_v^u \int_{v'}^u \partial_u \partial_v \phi(u',v') du' dv' & =
\phi(u,u) - \phi(u,v) - \int_v^u \partial_v \phi(v',v') dv'.
\label{easy}
\end{align}
Hence
\begin{align}
\Norm{\phi}_{L^{\infty}_{u,v}}  & \leq \Norm{\phi(x,0)}_{L^{\infty}_x} +
\Norm{\partial_t \phi(x,0)}_{L^1_x} + 
\Norm{\partial_x \phi(x,0)}_{L^{1}_x} + 
\Norm{\partial_u \partial_v \phi}_{L^1_{u,v}} \nonumber \\
&   \leq \Norm{\phi}_X.  \label{linfty.bound}
\end{align}

To prove \eqref{algebra-property}, it suffices to bound $\partial_u \partial_v (\phi \psi)$ in
$L^1_{u,v}$, since $\dd$ is easily seen to be an algebra.  We compute using
\eqref{mixed.est} and \eqref{linfty.bound}:
\begin{align*}
\Norm{\partial_u \partial_v(\phi \psi)}_{L^1_{u,v}} & \leq
\Norm{\partial_u \partial_v\phi \cdot  \psi}_{L^1_{u,v}}
+ \Norm{\partial_u\phi \cdot \partial_v\psi}_{L^1_{u,v}} \\
&  + \Norm{\partial_v \phi \cdot \partial_u \psi}_{L^1_{u,v}} +
\Norm{\partial_u \partial_v \psi \cdot \phi}_{L^1_{u,v}} \\
& \leq  \Norm{\partial_u \partial_v \phi}_{L^1_{u,v}}
\Norm{\psi}_{L^{\infty}_{u,v}} +
\Norm{\partial_u \phi}_{L^{1}_u L^{\infty}_v}
\Norm{\partial_v \psi}_{L^1_vL^{\infty}_u}  \\  
& + \Norm{\partial_v \phi}_{L^1_vL^{\infty}_u}
\Norm{\partial_u \psi}_{L^{1}_u L^{\infty}_v} + 
\Norm{\partial_u \partial_v \psi}_{L^1_{u,v}}
\Norm{\phi}_{L^{\infty}_{u,v}}    \\
& \lesssim \Norm{\phi}_X \Norm{\psi}_X 
\end{align*}

It remains to prove \eqref{box-property}.  We compute:
\begin{align*}
\Norm{\Box^{-1}\big( \phi \Box \psi \big)}_{X}  
& = \Norm{\partial_u \partial_v 
\big( \Box^{-1} ( \phi \Box \psi )\big)}_{L^1_{u,v}} \\
& = \Norm{ \phi \partial_u \partial_v
 \psi}_{L^1_{u,v}} \\
& \leq \Norm{\phi}_{L^{\infty}_{u,v}} 
\Norm{\partial_u \partial_v \psi}_{L^1_{u,v}} \\
& \leq \Norm{\phi}_{X} \Norm{\psi}_X.
\end{align*} 

This concludes the proof of global well-posedness for data 
which is small in $\dd$
and whose image lies in a single coordinate chart.  The following 
elementary lemma allows us to localize large data.  Let $\chi(x)$ be 
a smooth bump function supported on $[-2,2]$ with $\chi = 1$ 
on $[-1,1]$.    
\begin{lemma} \label{standard.local} Given
$\epsilon > 0$ and data $(f,g) \in L^{1,1}$
there exists $\delta = \delta(f,g,\epsilon)$ so that for
all $x_0 \in \rr$, $\chi(\frac{x - x_0}{\delta}) f $ takes values 
in a single coordinate chart of $\m$ and    
\begin{align}
\label{localizeI}
\|\chi(\frac{x - x_0}{\delta})(f,g)\|_{\dd} & \leq \epsilon.
\end{align}
\end{lemma}
\begin{proof} Fix $x_0 \in \rr$, and pick $\delta$ so that
\begin{align}
\int_{|x - x_0| \leq \delta} \norm{\partial_x f} dx +
\int_{|x - x_0| \leq \delta} \norm{g} dx & < \epsilon
\label{l.11.concentration}
\end{align}
for all $x_0 \in \rr$.  We may assume a coordinate chart 
on $\m$ around $f(x_0)$ is centered at $0$. Then together with
\eqref{l.11.concentration}, the fundamental theorem of calculus gives  
\begin{align*}
\|\chi(\frac{x - x_0}{\delta})f\|_{L^{\infty}_x} +
\|\chi(\frac{x - x_0}{\delta})g\|_{L^1_x}
&  \leq \epsilon. 
\end{align*} 
Finally,
\begin{align*}
\|\frac{d}{dx} \big( \chi(\frac{x - x_0}{\delta}) f\big)\|_
{L^1_x} & \lesssim
\| \frac{1}{\delta} f\|_{L^1(|x - x_0| \leq \delta)} +
\|f'\|_{L^1(|x - x_0| \leq \delta)} \\
& \leq \|f\|_{L^{\infty}(|x - x_0| \leq \delta)} + \epsilon \\
& \lesssim \epsilon
\end{align*}
\end{proof}

Together with the small data argument providing both 
existence and uniqueness of solutions, 
Lemma \ref{standard.local}  and finite speed of propagation
give local well-posedness for large $L^{1,1}$ data, but with a 
time of existence depending upon how concentrated the data is 
in $L^{1,1}$. 

\subsection{Nonconcentration of $L^{1,1}$ norm}

We now turn to the question of global well-posedness for large data.  
Suppose for contradiction that there exists large $L^{1,1}$ 
data $f$, $g$ for which a solution $\phi$
could only be continued in $X$ up to a maximal time of existence $0 < T^* < 
\infty$.  By finite speed of propagation we may assume that $\phi$ is compactly
supported.
 
From the small-data well-posedness theory (which in particular implies uniqueness), finite speed 
of propagation, and a Lemma \ref{standard.local} 
this implies the existence of a point $x_0$ such that the solution
concentrates on intervals near $x_0$ on every coordinate chart\footnote{Note
that we are using the hypothesis that the Christoffel symbols are uniformly
analytic to make these estimates independent of the choice of chart.}:
$$ \lim \sup_{\tau \to 0} \| \phi(T^*-\tau) \|_{\dd(x_0-4\tau, x_0+4\tau)} \geq \delta > 0.$$
By translation invariance we may take $x_0 = 0$.
We may pick our coordinate charts at each time $T^*-\tau$ so that $\phi(0, T^*-\tau) = 0$.
By the fundamental theorem of calculus, the concentration thus becomes
$$ \lim \sup_{\tau \to 0} \| \phi_x(T^*-\tau) \|_{L^1(-4\tau, 4\tau)}
+ \| \phi_t(T^*-\tau) \|_{L^1(-4\tau, 4\tau)} \gtrsim \delta > 0.$$
We can rewrite these derivatives in terms of $u$, $v$ derivatives to obtain
$$ \lim \sup_{\tau \to 0} \| \phi_u(T^*-\tau) \|_{L^1(-4\tau, 4\tau)}
+ \| \phi_v(T^*-\tau) \|_{L^1(-4\tau, 4\tau)} \gtrsim \delta > 0.$$
Thus in order to obtain a contradiction we need only show that the $L^1$ norms of
$|\phi_u|_\h$ and $|\phi_v|_\h$ do not concentrate.

But this is an immediate consequence of Lemma \ref{pohl_law}.
Indeed, as the data is in $L^{1,1}$,
the quantities $|\phi_u|_\h$ and $|\phi_v|_\h$ are travelling waves of $L^1$ functions
and therefore do not concentrate.

\subsection{Conformal compactification and scattering}
\label{l11.scatter}
Let $\phi$ denote a global $L^{1,1}$ solution to the wave map equation \eqref{wavemap}, and let
$\phi^+$, $\phi^-$ denote global $L^{1,1}$ solutions to the 
free wave equation.  Note that these solutions are continuous, 
since the solution space $X$ used earlier embeds into the space of 
continuous functions.  To show scattering and asymptotic completeness 
we have to prove the following two claims:
\begin{itemize}
\item For any $\phi$, there exists a $\phi^+$ such that $\| \phi(T) - \phi^+(T) \|_{L^{1,1}} \to 0$
as $T \to \infty$.
\item For any $\phi^-$, there exists a $\phi$ such that $\| \phi(T) - \phi^-(T) \|_{L^{1,1}} \to 0$
as $T \to -\infty$.
\end{itemize}
In $\rr^{1+1}$, the conformal compactification transformation (see \cite{christ.null, penrose})
takes the form
$$ (P\phi)(U,V) = \phi(\tan U, \tan V)$$
which takes functions on $\R^{1+1}$ to functions on the Einstein diamond 
$\{|U|, |V| < \pi/2\}$.
Since the wave map equation and the free wave equation are both conformally invariant in
one dimension, the function $P\phi$ is also a solution to \eqref{wavemap},
and  $P\phi^+$, $P\phi^-$ are solutions to the free wave equation.

A quick computation shows that the $L^{1,1}$ norm of $(P\phi)(0)$ is 
equal to the
$L^{1,1}$ norm of $\phi(0)$ since the Jacobian
factor and chain rule factor cancel.  Thus by the global 
well-posedness theory just proved,
$(P\phi)$ extends to a solution $\Phi$ to \eqref{wavemap} on all of $\R^{1+1}$,
where we may continuously extend the initial data so that the initial position $\Phi(0)$ is
constant and the initial velocity $\Phi_t(0)$ is zero on the intervals $(-\infty,-\pi/2)$
and $(\pi/2,\infty)$.

By Lemma \ref{scatter}, we
see that $\Phi$ is exactly equal to an $L^{1,1}$ 
solution $\Phi^+$
(resp. $\Phi^-$) to the free wave equation for $U \geq \pi/2$ or
$V \geq \pi/2$ (resp. $U \leq -\pi/2$ or $V \leq -\pi/2$). See
Figure \ref{fig:diamond}. 
We may of course extend $\psi^\pm$ to be $L^{1,1}$ solutions to the free wave 
equation on all of $\R^{1+1}$.  This gives a well-defined map from
$\Phi$ to $\Phi^\pm$; the corresponding inverse map also exists by the
same reasoning.

We now define the scattering maps $W_\pm: \phi \to \phi_\pm$ by defining
$\phi^\pm = L^{-1} \Phi^\pm$ on the Einstein diamond; it is easy
to see from the above discussion
that this map is well-defined and invertible, and that
 $\phi^\pm$ are global $L^{1,1}$ solutions to the free wave 
equation.  To complete the proof of scattering it suffices to show
that
$$\|\phi(T) - \phi^\pm(T)\|_{L^{1,1}} \to 0$$ as $T \to \pm \infty.$

By time reversal symmetry it suffices to do this for $\phi^+$.  Since by construction
$\phi^+(T)$ and $\phi(T)$ agree at the boundary of the Einstein diamond
(i.e. when $x \to \pm \infty$), 
it suffices by the fundamental theorem of calculus to show that
$$ \| (\phi - \phi^+)_u(T) \|_{L^1_x} + 
\| (\phi - \phi^+)_v(T) \|_{L^1_x} \to 0$$
as $T \to \infty$.  We show this only for the first term, as the second
is analogous.  We have
\begin{equation}
\label{twopage}
 \| (\phi - \phi^+)_u(T) \|_{L^1_x} = \int |(\phi - \phi^+)_u(u, u-2T)|\ du.
\end{equation}
By changing to the $U$ and $V$ coordinates, \eqref{twopage}
is  
$$\int_{-\pi/2}^{\pi/2} |(\Phi-\Phi^+)_U(U, 
\tan^{-1}(\tan(U) - 2T)|\ dU.$$
By the fundamental theorem of calculus and the fact that $\Phi-\Phi^+$ vanishes
at the upper boundary of the Einstein diamond, this is majorized by
$$\int\int_{|U|,|V| \leq \pi/2, V > \tan^{-1}(\tan(U)-2T)}
|(\Phi-\Phi^+)_{UV}(U,V)|\ dU dV.$$
By the monotone convergence theorem, this will go to zero as $T \to 0$ provided
that
$$\int\int_{|U|,|V| \leq \pi/2} |(\Phi-\Phi^+)_{UV}(U,V)|\ dU dV < \infty.$$
But since $\Phi$ obeys \eqref{wave-map}, $|\Phi_{UV}| \lesssim |\Phi_U| 
|\Phi_V|$.  Since $\Phi^+$ is a free solution, $\Phi^+_{UV} = 0$ and
so our integral is majorized by
$$\int\int_{|U|,|V| \leq \pi/2} |\Phi_U|_\h |\Phi_V|_\h\ dU dV,$$
But by Lemma \ref{pohl_law} this is majorized by
the square of the $L^{1,1}$ norm of the data of $\Phi$, which is finite.

\end{document}